\documentclass[11pt]{article}
\usepackage{amssymb}
\usepackage{amsmath}
\usepackage{mathrsfs}
\usepackage{graphics}
\usepackage{graphicx}
\usepackage{epstopdf}
\usepackage{extarrows}
\usepackage{xcolor}
\usepackage{subfigure}
\usepackage[T1]{fontenc}
\usepackage{latexsym,amssymb,amsmath,amsfonts,amsthm}
\usepackage{txfonts}
\usepackage{lineno}
\newcommand{\be}{\begin{eqnarray}}
\newcommand{\ben}{\begin{eqnarray*}}
\newcommand{\en}{\end{eqnarray}}
\newcommand{\enn}{\end{eqnarray*}}

\newtheorem{theorem}{Theorem}[section]
\newtheorem{lemma}{Lemma}[section]
\newtheorem{prp}[theorem]{Proposition}

\newtheorem{cor}[theorem]{Corollary}
\newtheorem{dfn}{Definition}[section]
\newtheorem{exmp}[theorem]{Example}
\newtheorem{remark}{Remark}

\newtheorem{Hypothesis}[theorem]{Hypothesis}

\definecolor{xty}{rgb}{1,0.5,0}
\begin{document}
\renewcommand{\theequation}{\arabic{section}.\arabic{equation}}
\begin{titlepage}
\title{\bf On Limit Measures and Their Supports for Stochastic Ordinary Differential Equations}
\author{\ \ Tianyuan Xu$^{1}$, Lifeng Chen$^{1}$, \ Jifa Jiang$^{1}$ $\footnote{Corresponding author.}$\\
{\small $^1$ Mathematics and Science College, Shanghai Normal University, Shanghai 200234, PR China.}\\
({\small {\sf XTY\_work@outlook.com}, \ {\sf lfchen@shnu.edu.cn},\ {\sf jiangjf@shnu.edu.cn},\ })}
\date{}
\end{titlepage}
\maketitle

\noindent\textbf{Abstract}:
 This paper studies limit measures of stationary measures of stochastic ordinary differential equations on the Euclidean space and tries to determine which invariant measures of an unperturbed system will survive. Under the assumption for SODEs to admit  the Freidlin--Wentzell or Dembo--Zeitouni large deviations principle with weaker compactness condition, we prove that limit measures are concentrated away from repellers  which are topologically  transitive, or equivalent classes, or admit Lebesgue measure zero.  We also preclude concentrations of limit measures on acyclic saddle or trap  chains. This illustrates that limit measures are concentrated on Liapunov stable compact invariant sets. Applications are made to the Morse--Smale systems, the Axiom A systems including structural stability systems and separated star systems, the gradient or gradient-like systems, those systems possessing the Poincar\'{e}--Bendixson property with a finite number of limit sets to obtain that limit measures live on Liapunov stable critical elements, Liapunov stable basic sets,  Liapunov stable equilibria and Liapunov stable limit sets including equilibria, limit cycles and saddle or trap cycles, respectively.  A number of nontrivial examples admitting a unique limit measure are provided, which include monostable, multistable systems and those possessing infinite equivalent classes.

\noindent \textbf{AMS Subject Classification}:\ \ Primary 60H10,37B35,60B10;
Secondary 60F10, 37A50, 37C70.

\noindent\textbf{Keywords}: large deviations; stationary measure; limit measure; concentration; Morse-Smale system; Axiom A system; gradient-like system; Poincar\'{e}-Bendixson property.

\section{Introduction}
This paper studies limit measures and their supports of stationary measures for stochastic ordinary differential equations
\begin{equation}\label{itodff}
  \mathrm dX^{\varepsilon}_t=b(X^{\varepsilon}_t)\mathrm dt+\varepsilon\sigma(X^{\varepsilon}_t)\mathrm dw_t, \quad X^{\varepsilon}_0=x\in \mathbb{R}^r
\end{equation}
when $\varepsilon$ goes to zero, where $w_t=(w^{1}_t,\cdots,w^{r}_t)^*$ is a standard $r$-dimensional Wiener process, the diffusion matrix $a=(a_{ij})_{r\times r}=\sigma\sigma^*$ is positive definite, which is called to be nondegenerate,
where $^{\ast}$ denotes transpose.
System (\ref{itodff}) is regarded as a stochastic perturbation of the deterministic dynamical system
\begin{equation}\label{unpersys}
  \frac{\mathrm dx}{\mathrm dt}=b(x), \quad x(0)=x\in \mathbb{R}^r.
\end{equation}
This asymptotic problem was first proposed by Kolmogorov in the 1950s (see \cite[p.838]{Sinai}). Khasminskii \cite{Khasminskii} proved that a stationary measure of a Markov process on a torus converges weakly to an invariant measure for a dynamical system on the torus as diffusion term tends to zero, hence he has been regarded as the first one to realize that limit measures of stationary measures are invariant measures for unperturbed systems. Using large deviations principle, Freidlin and Wentzell solved the problem of rough asymptotics of stationary measures of a diffusion process with small diffusion on a {\bf compact connected  manifold} \cite{FW1, FW2}. Precisely,  under assumptions that there are a finite number of equivalent classes containing all limit sets of the unperturbed system, they created a method to determine which stable equivalent class or classes a limit measure of stationary measures is concentrated on after a series of rare probability estimates. Kifer \cite{Kifer} generalized the corresponding results of Freidlin and Wentzell to discrete-time dynamical systems on a compact manifold with small unbounded random perturbations satisfying large deviations principle. Ruelle \cite{Ruelle} verified that limit measures of stationary measures for discrete-time dynamical systems with small bounded random perturbations live on quasiattractors. Bena\"{\i}m \cite{Benaim1} investigated the dynamical properties of a class of urn processes and recursive stochastic algorithms with constant gain and proved that limit measures of the process are concentrated on the Birkhoff center of irreducible attractors of its averaging ordinary differential equations.

As for stochastic system (\ref{itodff}) on the Euclidean space, Freidlin and Wentzell \cite{FW2} and Hwang \cite{Hwang1980} considered gradient systems (that is, $b(x)=-\nabla U(x)$) perturbed by an additive noise and showed that limit measures of stationary measures have their supports on the lowest energy points via the large deviation technique and Laplace's method, respectively. Huang et al. \cite{Huang2018} proved that all limit measures of stationary measures of (\ref{itodff}) are invariant with respect to the solution flow of (\ref{unpersys}) and sit on the global attractor of (\ref{unpersys}) by estimating measure values of regular stationary measures in an exterior domain with respect to diffusion and Liapunov-like functions.
For a given deterministic ODEs (\ref{unpersys}) with a strongly local attractor (repeller), they constructed a nondegenerate diffusion $\sigma$  such that  limit measures of (\ref{itodff}) are concentrated on the local attractor (away from the local repeller); then Ji et al. \cite{Ji2019}  removed the ``strong'' hypothesis of attractor and repeller and obtained the same conclusion. Besides, Huang et al. \cite{Huang2018} proved all limit measures of (\ref{itodff}) for all small nondegenerate diffusion perturbations are concentrated away from  any hyperbolic repelling equilibrium by constructing a positive definite quadratic form as an anti-Liapunov function. Chen, Dong and Jiang \cite{Chen2020} presented a criterion that limit measures are concentrated away from repellers, which  can be applied to repelling limit cycles or quasi-periodic orbits. For a gradient system  (\ref{unpersys}),  Huang et al. \cite{Huang2016} proved all limit measures of (\ref{itodff}) for all small nondegenerate diffusion perturbations support on the set of critical points of the potential function $U$. Chen et al. \cite{ChenXinfu} analyzed limit measures of one dimensional system (\ref{itodff}) as the white noise vanishes and proved that all limit measures are exactly concentrated on the global minimizers of $-\int \frac{b(u)}{\sigma^2(u)}\mathrm du$ if $b'(x)\neq 0$ at all these global minimizers. This shows that the limit measure may support on those lesser stable equilibria but with smaller $\sigma$.

Under the assumptions that stochastic system (\ref{itodff}) admits  the Freidlin--Wentzell or Dembo--Zeitouni large deviations principle with weaker compactness condition, we shall exploit limit measures of stationary measures of stochastic ordinary differential equations (\ref{itodff}). Such  measures are more stable than other invariant measures of unperturbed systems (\ref{unpersys}) or the most stable if they uniquely exist to stochastic perturbations.  We shall prove that limit measures are concentrated away from repellers  which are topologically  transitive, or equivalent classes, or admit Lebesgue measure zero.  We also preclude concentrations of limit measures on acyclic saddle or trap  chains. This shows that limit measures are concentrated on Liapunov stable compact invariant sets. Applications are made to the Morse--Smale systems, the Axiom A systems including structural stability systems and separated star systems, the gradient or gradient-like systems, those systems possessing the Poincar\'{e}--Bendixson property with a finite number of limit sets to obtain that limit measures live on Liapunov stable critical elements, Liapunov stable basic sets,  Liapunov stable equilibria and Liapunov stable limit sets including equilibria, limit cycle and saddle or trap cycles, respectively.
As far as we know, there are seldom examples of SODEs (\ref{itodff}) on the Euclidean space whose limiting measures
and their supports are clearly described. In Section 5, a number of nontrivial examples admitting
a unique limit measure are provided, which include monostable, multistable systems
and those possessing infinite equivalent classes.

\section{Preliminaries and Notations}

In this section, we recall some basic definitions and preliminary results from
\cite{FW1, FW2, Kifer, DZ, Feng, Ruelle, Palis, Wen, Hartman}.
Throughout this paper we always assume that the coefficients  $b$ and $\sigma$ are locally Lipschitz continuous on $\mathbb{R}^r$, the solutions of (\ref{itodff}) and (\ref{unpersys}) are defined on $[0, +\infty)$, and the diffusion matrix $a=\sigma\sigma^*$ is nondegenerate on $\mathbb{R}^r$.

For $x\in\mathbb{R}^r,\delta>0$ and $\mathcal{S},\mathcal{T}\subset\mathbb{R}^r$,
we may write
$B(x,\delta)=B_{\delta}(x)=\{y\in\mathbb{R}^r:|y-x|<\delta\}$,
$\bar{B}_{\delta}(x)
=\{y\in\mathbb{R}^r:|y-x|\leq\delta\}$, ${\rm dist}(\mathcal{S},\mathcal{T})=\inf_{y\in\mathcal{S}, z\in\mathcal{T}}|y-z|$
and $(\mathcal{S})_{\delta}=\{y\in\mathbb{R}^r:{\rm dist}(y,\mathcal{S})<\delta\}$.
Here $|\cdot|$ denotes the usual Euclidean norm.

We recall some standard definitions and results of dynamical systems (see, e.g., \cite{Ruelle, Benaim1, Palis, Wen, Hartman}).
The solution semiflow of deterministic system (\ref{unpersys}) is denoted by $\Psi_{t}(x)=\Psi(t,x)$,
its positive (resp. negative, entire) orbit is denoted by $\gamma^+(x)$ (resp. $\gamma^-(x), \gamma(x)$),
and $\omega$-limit set (resp. $\alpha$-limit set) is denoted by $\omega(x)$ (resp. $\alpha(x)$).

For $D\subset\mathbb{R}^r$ and $T\geq 0$, let $\gamma^+(D)$ denote the sum of all positive orbits passing through points in $D$,
the notations $\Psi_{T}(D)$ and $\Psi([0,T]\times D)$ are defined in a similar manner.
A set $\Lambda\subset\mathbb{R}^r$ is {\it positively invariant} if $\Psi_{t}(\Lambda)\subset \Lambda$ for all $t\geq 0$.
It is {\it invariant} if $\Psi_{t}(\Lambda)=\Lambda$ for all $t\geq 0$.

A set $\Lambda\subset\mathbb{R}^r$
is called \emph{topologically transitive} (resp. \emph{minimal}) for $\Psi$ if $\Lambda$ is nonempty
compact invariant, and $\exists x\in\Lambda$ (resp. $\forall x\in\Lambda$)
such that $\overline{\gamma^{+}(x)}=\Lambda$. Obviously, any minimal set
is topologically transitive from the definitions.
Suppose that $\Lambda$ is topologically transitive. Then it follows
from the definition that for all open sets $U_1$ and $U_2$
in $\mathbb{R}^r$ such that $U_1\cap \Lambda$ and $U_2\cap \Lambda$ are nonempty,
there is $T>0$ such that $\Psi_{T}(U_1)\cap U_2\neq\emptyset$.

A subset $\mathcal{R}\subset \mathbb{R}^r$ is called to be a {\it repeller} (an {\it attractor}) for $\Psi$ provided:
(i) $\mathcal{R}$ is nonempty, compact and invariant; and
(ii) $\mathcal{R}$ has a neighborhood $N\subset \mathbb{R}^r$, called a {\it fundamental neighborhood} of $\mathcal{R}$,  such that
$\lim_{t\rightarrow-\infty}{\rm dist}\big(\Psi_{t}(x),\mathcal{R}\big)=0$ ($\lim_{t\rightarrow +\infty}
{\rm dist}\big(\Psi_{t}(x),\mathcal{R}\big)=0$) uniformly in $x\in N$. In the case of attractor, we call $\mathcal{R}$ to be
{\it Liapunov stable}.

\begin{prp}\label{reppp}
Let $\mathcal{R}$ be a repeller for $\Psi$ with a fundamental neighborhood $N$.
Then for any
compact set $K\subset N\backslash \mathcal{R}$,
there exists $T>0$ such that $\Psi_{t}(K)\subset N^c$ for all $t\geq T$.
\end{prp}
\begin{proof}
Let $K\subset N\backslash \mathcal{R}$ be a compact set.
Then $d:={\rm dist}(K,\mathcal{R})>0$.
Choose $\eta\in(0,d)$. Then $(\mathcal{R})_{\eta}\subset K^c$.
Since $\mathcal{R}$ is a repeller,
there is a $T>0$ such that for any $t\geq T$,
$\Psi_{-t}(N)\subset (\mathcal{R})_{\eta}\subset K^c$.
This implies that $N\subset \Psi_{t}(K^c)$ for any $t\geq T$.
Therefore, $\Psi_{t}(K)\subset N^c$ for any $t\geq T$.
\end{proof}

The dynamical system $\Psi$ on $\mathbb{R}^r$ is called \emph{dissipative} if there exists a bounded set $B\subset\mathbb{R}^r$
with the property that for every compact $K\subset\mathbb{R}^r$ there exists $T=T(K)>0$ such
that $\Psi_{t}(K)\subset B$ for all $t\geq T$. It is well known that the above dissipativity of $\Psi$
is equivalent to that $\Psi$ has a global attractor $\mathcal{A}$, that is,
$\mathcal{A}$ is an attractor whose basin is all the space $\mathbb{R}^r$.

We now review some notions and properties in the large deviation theory, which are useful
in deal with stationary measure asymptotics for small noise Markov processes
(see, e.g., \cite{FW1, FW2, Kifer}).

For each fixed $T>0$,
let ${\bf C}_{T}=C([0,T],\mathbb{R}^r)$ (resp. ${\bf AC}_T=AC([0,T],\mathbb{R}^r)$)
denote the set of continuous functions (resp. absolutely continuous functions) on $[0,T]$ with values in $\mathbb{R}^r$.
For $\varphi,\psi\in{\bf C}_{T}$ and $W\subset {\bf C}_{T}$, let
$\rho_{T}(\varphi,\psi)=\sup_{0\leq t\leq T}|\varphi(t)-\psi(t)|$
and $\rho_{T}(\varphi,W)=\inf_{\phi\in W}\rho_{T}(\varphi,\phi)$.
We define the following functional on ${\bf C}_{T}$:
\begin{equation}\label{actfal}
  S_{T}(\varphi)=S_{0T}(\varphi)
=\left\{
\begin{array}{ll}
\frac{1}{2}\int\limits_{0}^{T}\Big(\dot{\varphi}(t)-b(\varphi(t))\Big)^{*}a^{-1}\big(\varphi(t)\big)
\Big(\dot{\varphi}(t)-b(\varphi(t))\Big) {\rm d}t, & \hbox{if $\varphi\in{\bf AC}_T$,} \\
+\infty, & \hbox{otherwise.}
\end{array}
\right.
\end{equation}

Let $X^{\varepsilon,x}_{\cdot}$, parametrized by $x\in\mathbb{R}^r$ and $\varepsilon>0$,
be the solution of SDEs {\rm (\ref{itodff})}. For each $T>0$, $\varepsilon>0$ and $x\in\mathbb{R}^r$,
one can regard $X^{\varepsilon,x}$ as a ${\bf C}_{T}$-valued random variable.
Furthermore, to emphasize the dependence of initial conditions $x$, we also introduce the following notions
$\ {\bf C}_{T}^{x}=\{\varphi\in{\bf C}_{T}:\varphi_0=x\},
{\bf AC}_T^x=
\{\varphi\in{\bf AC}_T:\varphi_0=x\}$, and
functional $S_{T}^{x}$ on ${\bf C}_{T}$ :
\begin{equation}\label{actfalini}
  S_{T}^{x}(\varphi)=S_{0T}^{x}(\varphi)
=\left\{
\begin{array}{ll}
\frac{1}{2}\int\limits_{0}^{T}\Big(\dot{\varphi}(t)-b(\varphi(t))\Big)^{*}a^{-1}\big(\varphi(t)\big)
\Big(\dot{\varphi}(t)-b(\varphi(t))\Big) {\rm d}t, & \hbox{if $\varphi\in{\bf AC}_T^x$,} \\
+\infty, & \hbox{otherwise.}
\end{array}
\right.
\end{equation}
By the definition of $S_{T}^{x}$ in (\ref{actfalini}) we easily have the following assertion.
\begin{remark}\label{Rem1}
$S_{0T}^{x}(\varphi)=0$ if and only if $\varphi$ {\rm(}up to time $T${\rm)}
coincides with the solution $\Psi_{\cdot}(x)$
of the deterministic system {\rm(\ref{unpersys})}.
\end{remark}

We give several definitions of uniform large deviations principles
that are found in the literature.
Let $\mathcal{K}$ be a collection of all compact subsets of $\mathbb{R}^r$
and $\mathbb{F}_{T}^{x}(s)
=\{\varphi\in{\bf C}_{T}^x:S_{0T}^{x}(\varphi)\leq s \}$ for $T>0,x\in\mathbb{R}^r,s\geq 0$.
The first definition of a uniform large deviations principle
presented here is due to
Freidlin--Wenztell \cite{FW1} (see also \cite[p.74]{FW2}).
\begin{dfn}[Freidlin--Wentzell uniform large deviations principle over $\mathcal{K}$]
For any fixed $T>0$,
we say that the random variables $\{X^{\varepsilon,x}\}$
satisfy a Freidlin--Wentzell uniform large deviations principle with respect to
the functionals $S_{T}^{x}$  uniformly over $\mathcal{K}$, if

${\bf(I_u)}$ for each $s_0>0,\ \delta>0,\ \gamma>0$ and $K\in\mathcal{K}$ there exists $\varepsilon_0>0$
such that
\begin{equation}\label{ufldplbb}
  \mathbb{P}_{x}\{\rho_{T}(X^{\varepsilon},\varphi)<\delta\}\geq \exp\{-\varepsilon^{-2}\big(S_{T}(\varphi)+\gamma\big) \}
\end{equation}
for all $\varepsilon\in(0,\varepsilon_0],\ x\in K$ and $\varphi\in\mathbb{F}_{T}^{x}(s_0)$;

${\bf(II_u)}$ for each $s_0>0,\ \delta>0,\ \gamma>0$ and $K\in\mathcal{K}$ there exists $\varepsilon_0>0$
such that
\begin{equation}\label{ufldpupbb}
  \mathbb{P}_{x}\{\rho_{T}(X^{\varepsilon},\mathbb{F}_{T}^{x}(s))\geq\delta\}\leq \exp\{-\varepsilon^{-2}(s-\gamma)\}
\end{equation}
for all $\varepsilon\in(0,\varepsilon_0],\ s\leq s_0$ and $x\in K$.
\end{dfn}

For any $F\subset {\bf C}_{T}$,
let $S_{T}^{x}(F):=\inf_{\varphi\in F}S^{x}_{0T}(\varphi)$.
The next definition of uniformly large deviations principle
is given in Dembo--Zeitouni \cite[p.~216, Corollary 5.6.15]{DZ}).
\begin{dfn}[Dembo--Zeitouni uniform large deviations principle over $\mathcal{K}$]
For any $T>0$, we say that the random variables $\{X^{\varepsilon,x}\}$
satisfy a Dembo--Zeitouni uniform large deviations principle with respect to
the functionals $S_{T}^{x}$  uniformly over $\mathcal{K}$, if

${\bf(I_u')}$ for any $K\in\mathcal{K}$ and open $G\subset{\bf C}_{T} $,
\begin{equation}\label{ufldplbb2}
  \liminf_{\varepsilon\rightarrow 0}
\inf_{x\in K}\varepsilon^2\log\mathbb{P}_{x}(X^{\varepsilon}\in G)\geq -\sup_{x\in K}S_{T}^{x}(G),
\end{equation}
which implies that for any $K\in\mathcal{K}$, open $G\subset{\bf C}_{T} $ and $\gamma>0$,
there exists $\varepsilon_0>0$ such that for any $\varepsilon\in(0,\varepsilon_0], x\in K$
\begin{equation}\label{ufldplbb3}
\mathbb{P}_{x}(X^{\varepsilon}\in G)\geq\exp\left\{-\frac{\sup_{x\in K}S_{T}^{x}(G)+\gamma}{\varepsilon^2}\right\};
\end{equation}

${\bf(II_u')}$ for any $K\in\mathcal{K}$ and closed $F\subset{\bf C}_{T} $,
\begin{equation}\label{ufldpupbb2}
  \limsup_{\varepsilon\rightarrow 0}
\sup_{x\in K}\varepsilon^2\log\mathbb{P}_{x}(X^{\varepsilon}\in F)\leq -\inf_{x\in K}S_{T}^{x}(F),
\end{equation}
which implies that for any $K\in\mathcal{K}$, closed $F\subset{\bf C}_{T} $ and $\gamma>0$,
there exists $\varepsilon_0>0$ such that for any $\varepsilon\in(0,\varepsilon_0], x\in K$
\begin{equation}\label{ufldpupbb3}
 \mathbb{P}_{x}(X^{\varepsilon}\in F)\leq\exp\left\{-\frac{\inf_{x\in K}S_{T}^{x}(F)-\gamma}{\varepsilon^2}\right\}.
\end{equation}
\end{dfn}

\begin{remark}
Below we will often write {\rm FWULDP} and {\rm DZULDP} as shorthand for
Freidlin--Wentzell and Dembo--Zeitouni uniform large deviations principle respectively.
Since bounded sets in $\mathbb{R}^r$ have compact closure, we could also rewrite
uniformly over compact sets by bounded sets.
We refer to the recent paper {\rm \cite[Theorem 2.7]{Sals}} for more details concerning equivalence between {\rm FWULDP} and {\rm DZULDP}
under additional assumptions.
\end{remark}

Throughout the rest of the paper we assume that $S_T$ enjoys the following property.
\begin{Hypothesis}\label{closeinfH}
For any $T>0$, let $F\subset{\bf C}_{T} $ be a closed  set and $F^0:=\{\varphi(0):\varphi\in F \}$ a bounded set in $\mathbb{R}^r$. Then
\[  \inf_{\varphi\in F}S_T(\varphi)=S_T(F)>0\]
if $F$ does not contain any solution of system {\rm (\ref{unpersys})}.
\end{Hypothesis}

Now we introduce the so--called conditions ${\bf(0_c)}$ and ${\bf(0_w)}$ as follows.

${\bf(0_c)}$ $S_{T}$ is lower semi-continuous and
the set $\cup_{x\in K}\mathbb{F}_{T}^{x}(s)=\{\varphi\in{\bf C}_{T}:\varphi(0)\in K,S_{T}(\varphi)\leq s \}$
is compact for each $s<+\infty$ and $K\in\mathcal{K}$.

${\bf(0_w)}$ $S_{T}$ is lower semi-continuous and
the set $\{\varphi\in{\bf C}_{T}:\varphi(t)\in K,t\in [0,T],S_{T}(\varphi)\leq s \}$
is compact for each $s<+\infty$ and $K\in\mathcal{K}$.

It is easy to see that ${\bf(0_c)}$ implies ${\bf(0_w)}$.
\begin{prp}\label {Hypothesis} Hypothesis \ref{closeinfH} holds if either ${\bf(0_c)}$ or ${\bf(0_w)}$ with dissipativity of {\rm (\ref{unpersys})} is satisfied.
\end{prp}
\begin{proof} The proof of the first part can be found in \cite{FW2, DZ};
The proof of the second part is postponed to Appendix.
\end{proof}

In the subsequent contents, we always assume that the solution of system (\ref{itodff}) admits FWULDP or DZULDP,
we leave the conditions for FWULDP and DZULDP to hold open, readers can refer to \cite{FW2, DZ, Feng, Kraaija, Yang}
and many references therein.

\emph{Quasipotential}, introduced by Freidlin and Wentzell (see, e.g., \cite [p.90]{FW2}),
is a very useful notion and is defined by
\[  V(x,y):=\inf\big\{S_{T}(\varphi):\varphi(0)=x,\ \varphi(T)=y,\ T\geq 0\big\},\ \quad x,\ y\in\mathbb{R}^r.  \]
Without ambiguity we can also define $V$ on pairs of subsets of $\mathbb{R}^r$
\[  V(D_1,D_2):=\inf\big\{S_{T}(\varphi):\varphi(0)\in D_1,\varphi(T)\in D_2, T\geq 0\big\}
=\inf_{x\in D_1,y\in D_2}V(x,y), \ \quad D_1,D_2\subset\mathbb{R}^r.  \]
A set $K$ is called to be an {\it equivalent class} if $V(x,y)=0$ for any $x,y\in K$.
We note that every limit set is an equivalent class.

In the following,  we introduce the Linear Interpolation Function of $x,y\in\mathbb{R}^r$ (abbreviated as
LIF$_{xy}$):
\[  \textrm{LIF$_{xy}$}(t)= x+\frac{t}{|y-x|}(y-x), \: t\in[0,|y-x|].  \]

The following basic property of $V$ is useful in our paper.
\begin{lemma}\label{Vctin}
For each compact $K\subset \mathbb{R}^r$, there is a positive constant $L=L_K$
such that for any $x,y\in K$ there exists a $C^\infty$ function $\varphi$, $\varphi(0)=x,\varphi(|x-y|)=y$ for which
$S_{|x-y|}(\varphi)\leq L|x-y|$.
In particular, $V(x,y)\leq L|x-y|$.
\end{lemma}
\begin{proof}
For instance, to choose $\varphi$ as LIF$_{xy}$ suffices.
\end{proof}

\begin{dfn}
Let $T_i>0,\ \varphi_i\in {\bf C}_{T_i},\ i=1,2$ and $\varphi_1(T_1)=\varphi_2(0)$.
We define a link function of $\varphi_1$ and $\varphi_2$, $\varphi_1\ast \varphi_2:[0,T_1+T_2]\longrightarrow \mathbb{R}^r$ by
\[  \varphi_1\ast \varphi_2(t)=
\left\{
  \begin{array}{ll}
    \varphi_1(t), & t\in[0,T_1]; \\
    \varphi_2(t-T_1), & t\in[T_1,T_1+T_2].
  \end{array}
\right.
  \]
The link function of a finite number of functions is defined similarly.
\end{dfn}

\section{Main Results}
First, we introduce a condition on a set $\mathcal{R}$ as follows.

${\bf(P_{\mathcal{R}})}$ For every $\eta>0$, there exist an open neighborhood $U$ of $\mathcal{R}$
and a positive constant $T^*$ such that for each $x\in\mathcal{R}$,  there exist $T=T(x)\leq T^*$ and $\varphi\in {\bf C}_{T}$ with $\varphi(0)=x$, $\varphi(T)\in U^c$ and
$S_{T}^{x}(\varphi)<\eta.$

This property depends only on the structure of the dynamical system (\ref{unpersys}). It will be proved that ${\bf(P_{\mathcal{R}})}$ holds for $\mathcal{R}$ either to have zero Lebesgue measure or to be topologically transitive or an equivalent class with a mild requirement.
\begin{theorem}\label{repeller}
Let the diffusion matrix $a$ be nonsingular on $\mathbb{R}^r$ and system {\rm(\ref{itodff})} satisfy Hypothesis \ref{closeinfH}. Suppose system {\rm(\ref{itodff})} possesses {\rm DZULDP} or {\rm FWULDP}.  If $\mathcal{R}$ is a repeller of {\rm(\ref{unpersys})} admitting the property ${\bf(P_{\mathcal{R}})}$
and $\mu^{\varepsilon}$  is any stationary distributions of diffusion process
$(X^{\varepsilon},\mathbb{P}_x)$,
then there exist a neighborhood $U_0$ of $\mathcal{R}$, $\kappa>0$ and $\varepsilon_*>0$
such that for any $\varepsilon\in(0,\varepsilon_*)$, we have
\[ \mu^{\varepsilon}(U_0)\leq \exp\{-\kappa/ \varepsilon^2\}  .   \]
As a result, if $\mu^{\varepsilon_j}\xlongrightarrow{w} \mu$ as $\varepsilon_j\rightarrow 0$,
then $\mu(U_0)=0$.
\end{theorem}

\begin{theorem}\label{saddle}
Let the diffusion matrix $a$ be nonsingular on $\mathbb{R}^r$ and system {\rm(\ref{itodff})} satisfy Hypothesis \ref{closeinfH} and possess {\rm DZULDP} or {\rm FWULDP}. Suppose that there are points $x_1, x_2, \cdots, x_m$ in $\mathbb{R}^r$, an attractor $\mathcal{A}=:\mathcal{S}_{m+1}$ and compact equivalent classes $\mathcal{S}_1, \mathcal{S}_2, \cdots, \mathcal{S}_m$  such that $\mathcal{A}$ and $\mathcal{S}_i\supset\alpha(x_i), i=1, 2,\cdots, m$ are pairwise disjoint, $\omega(x_i)\subset \mathcal{S}_{i+1}, i=1, 2,\cdots, m$.
 If
 $\mu^{\varepsilon}$ is any stationary distributions of diffusion process
$(X^{\varepsilon},\mathbb{P}_x)$, then there exist neighborhoods $U_i$ of $\mathcal{S}_i, i=1, 2, ..., m$, $\kappa>0$ and $\varepsilon_*>0$
such that for any $\varepsilon\in(0,\varepsilon_*)$, we have
\[ \mu^{\varepsilon}(U_i)\leq \exp\{-\kappa/ \varepsilon^2\}  .   \]
As a result, if $\mu^{\varepsilon_j}\xlongrightarrow{w} \mu$ as $\varepsilon_j\rightarrow 0$,
then $\mu(U_i)=0, i=1, 2, ..., m$.
\end{theorem}
Theorems \ref{repeller} and \ref{saddle} are proved in Section 4 below.

Note that Huang et al. \cite{Huang2018} and Ji et al. \cite{Ji2019} showed that for a given deterministic system {\rm(\ref{unpersys})} with a local repeller there exists a nondegenerate perturbed system {\rm(\ref{itodff})}   such that all its limit measures are concentrated away from the local repeller and only  Huang et al. \cite{Huang2018} got rid of the concentration on  any hyperbolic repelling equilibrium for all small nondegenerate diffusion perturbations. These results depend on diffusions $\sigma$ except hyperbolic repelling equilibrium. Theorem \ref{repeller} asserts that all limit measures for {\rm(\ref{itodff})} are concentrated away from any local repeller with the ${\bf\big(P_{\mathcal{R}}\big)}$, which is independent of nondegenerate diffusions. Our result is easier to use because the existing conditions for FWULDP is very weak and there is no need to construct Liapunov function for a given compact invariant set, which is not an easy job. Finally, we note that if ${\bf\big(P_{\mathcal{R}}\big)}$ is violated then Theorem \ref{repeller} does not hold, see Example \ref{counterexample}. This means the condition ${\bf\big(P_{\mathcal{R}}\big)}$ is sharp.

Theorem \ref{saddle} presents a criterion on noconcentration on saddles or semistable recurrent orbits or a continuum of stable recurrent orbits.  Chen et al. \cite{Chen2020, Chen2, Chen3} constructed examples concentrated on a saddle or saddles. Observing the examples, we find that the corresponding deterministic systems possess saddle cycle which is homoclinic orbit or saddle-connections forming a cycle. Theorem \ref{saddle} shows that the necessity  for a limit measure to concentrate on a saddle or saddles is that the  deterministic system {\rm(\ref{unpersys})} admits at least a saddle cycle. This precludes  concentration on saddles of a series of examples in \cite{Chen2} and solves the conjecture proposed in that paper.

Finally, we remark that as far as we know there are seldom examples of SODEs (\ref{itodff})
on the Euclidean space whose limiting measures and their supports are clearly described.
In Section 5, we will present a plenty of examples such that corresponding limiting measures
and their supports are clearly given by Theorems \ref{repeller} and \ref{saddle}.

\section{The Proofs of the Main Results}

We first prove two key lemmas which say the quasipotential against the flow of {\rm(\ref{unpersys})} is positive near repeller or attractor. Both lemmas play important roles in this paper.
\begin{lemma}\label{klofR}
Let the diffusion matrix $a$ be nonsingular on $\mathbb{R}^r$ and system {\rm(\ref{itodff})} satisfy Hypothesis \ref{closeinfH}. If $\mathcal {R}$ is a repeller of the system {\rm(\ref{unpersys})},
then for any $\delta>0$, there exist $\delta_1,\delta_2\in(0,\delta), \delta_1>\delta_2$
and $s_0>0$ such that $V(\partial (\mathcal{R})_{\delta_1},\partial(\mathcal{ R})_{\delta_2})\geq s_0$.
\end{lemma}
\begin{proof}
Since $\mathcal{R}$ is a repeller, we may take $\delta_1^*\in(0,\delta)$ such
that $(\mathcal{R})_{\delta_1^*}$ is a fundamental neighborhood of $\mathcal{R}$.
Let $\delta_2^*=\delta_1^* /2$. Since $\partial (\mathcal {R})_{\delta_2^*}\subset (\mathcal{R})_{\delta_1^*}\backslash \mathcal{R}$
is compact, by Proposition \ref{reppp}, there exists $T_0>0$ such that
\begin{equation}\label{useRep}
\Psi_{t}\left(\partial (\mathcal{R})_{\delta_2^*}\right)
\subset ((\mathcal{R})_{\delta_1^*})^c, \forall t\geq T_0.
\end{equation}
Define $d_0={\rm dist}\left(\Psi\big([0,T_0]\times \partial(\mathcal{R})_{\delta_2^*} \big),\mathcal{R}\right)$. Then by the invariance of $\mathcal{R}$ and the definition of $d_0$, we have
\begin{equation}\label{d0}
0<d_0\leq\delta_2^*.
\end{equation}
Let $\delta_1=\frac{3}{4}\delta_1^*, \delta_2=\frac{d_0}{2}$
and $\tilde{\delta}=\frac{d_0}{4}$.
Then $\tilde{\delta}<\delta_2\leq\frac{1}{4}\delta_1^*<\delta_2^*<\delta_1<\delta_1^*$ by (\ref{d0}).
$s_0$ can be taken as
\[ s_0:= \inf\left\{S_{T_0}(\varphi):\varphi\in {\bf C}_{T_0},\varphi(0)\in \partial (\mathcal{R})_{\delta_2^*},
\rho_{T_0}\big(\varphi(\cdot),\Psi_{\cdot}(\varphi(0))\big)\geq \tilde{\delta} \right\}.  \]
By Hypothesis \ref{closeinfH}, we get $s_0>0$. We claim that
$$V(\partial(\mathcal{R})_{\delta_1},\partial(\mathcal{R})_{\delta_2})\geq s_0.$$
In fact, if it is not true,
then there exist $T\in(0,+\infty)$ and $\psi\in {\bf C}_T, \psi(0)\in\partial(\mathcal{R})_{\delta_1},
\psi(T)\in\partial(\mathcal{R})_{\delta_2}$ such that
$S_T(\psi)<s_0$. Let $t^*:=\sup\{t\in[0,T]:\psi(t)\in \partial(\mathcal{R})_{\delta_2^*}\}$. Then
by the continuity of $\psi$
and $\delta_2<\delta_2^*<\delta_1$,
we get that $0<t^*<T$, $\psi(t^*)\in \partial(\mathcal{R})_{\delta_2^*}$ and $\psi(t)\notin \partial(\mathcal{R})_{\delta_2^*}, t\in (t^*, T].$ We extend the definition
of function $\psi$ up to the end of the interval $[0,t^*+T_0]$
as the solution of deterministic system {\rm(\ref{unpersys})} as long as $T < t^*+T_0$.  Thus $\psi(t^*+\cdot)|_{[0,T_0]}\in{\bf C}_{T_0}$. By the nonnegativity and additivity of $S$,
we have $S_{t^*,t^*+T_0}(\psi)\leq S_{t^*,T}\leq S_{T}(\psi)<s_0$.
From the definition of $s_0$, we get that
\begin{equation}\label{uses0}
\rho_{T_0}\big(\psi(t^*+\cdot),\Psi_{\cdot}(\psi(t^*))\big)< \tilde{\delta}.
\end{equation}

Case 1. $T_0\leq T-t^*$.
It follows from (\ref{useRep}) and (\ref{uses0}) that
\begin{align*}
  {\rm dist}\big(\psi(t^*+T_0),\mathcal{R}\big)\geq & -\left|\psi(t^*+T_0)-\Psi_{T_0}(\psi(t^*))\right|+{\rm dist}\big(\Psi_{T_0}(\psi(t^*)),\mathcal{R}\big)\\
  \geq &-\rho_{T_0}(\psi(t^*+\cdot),\Psi_{\cdot}(\psi(t^*)))
+{\rm dist}\big(\Psi_{T_0}(\psi(t^*)),\mathcal{R}\big)\\
  \geq & \delta_1^*-\tilde{\delta}\geq \frac{7}{8}\delta_1^*\\
  > &\delta_1.
  \end{align*}
The continuity of $\psi$ implies
that there is a time $\hat{t}\in(t^*+T_0,T)$
such that $\psi(\hat{t})\in \partial (\mathcal{R})_{\delta_2^*}$.
This contradicts the definition of $t^*$.

Case 2. $T_0> T-t^*$. For any $t\in[0,T_0]$, by (\ref{uses0}) and the definition of $d_0$, we have
\begin{align*}
  {\rm dist}\big(\psi(t^*+t),\mathcal{R}\big)\geq & -|\psi(t^*+t)-\Psi_{t}(\psi(t^*))|
  +{\rm dist}\big(\Psi_{t}(\psi(t^*)),\mathcal{R}\big)\\
  >&-\tilde{\delta}+d_0=\frac{3}{4}d_0\\
  >& \delta_2.
\end{align*}
 This contradicts to the
condition $\psi(t^*+T-t^*)=\psi(T)\in \partial(\mathcal{R})_{\delta_2}$ while $T-t^*\in(0,T_0)$.

This proves the claim and completes the proof.
\end{proof}

Similarly, we can prove the same version about an attractor.
\begin{lemma}\label{klofA}
Let the diffusion matrix $a$ be nonsingular on $\mathbb{R}^r$ and system {\rm(\ref{itodff})} satisfy Hypothesis \ref{closeinfH}. If $\mathcal{A}$
is an attractor of the system {\rm(\ref{unpersys})},
then for any $\delta>0$, there exist $\delta_1,\delta_2\in(0,\delta), \delta_1>\delta_2$
and $s_0>0$ such that $V(\partial (\mathcal{A})_{\delta_2},\partial (\mathcal{A})_{\delta_1})\geq s_0$.
\end{lemma}
\begin{proof}
First choose $\delta_1,\delta_1^*: \delta>\delta_1>\delta_1^*>0$ so small that $(\mathcal{A})_{\delta_1}$ is a fundamental neighborhood of $\mathcal{A}$. The definition of attractor implies that there is $\delta^*_2<\delta^*_1$ such that $\gamma^+\Big(\overline{(\mathcal{A})_{\delta^*_2}}\Big)\subset (\mathcal{A})_{\delta^*_1}$. Next choose $\delta_2,\delta^*_3: \delta^*_2>\delta_2>\delta_3^*>0$. As $(\mathcal{A})_{\delta_1}$ is a fundamental neighborhood of $\mathcal{A}$ and $\partial(\mathcal{A})_{\delta^*_2}\subset (\mathcal{A})_{\delta_1}$, there exists a $T_0>0$ satisfying $\gamma^+\big(\Psi_{T_0}(\partial(\mathcal{A})_{\delta^*_2})\big)\subset (\mathcal{A})_{\delta^*_3}$. From now on  we fix $\delta>\delta_1>\delta_1^*>\delta_2^*>\delta_2>\delta_3^*>0$ and set $\tilde{\delta}=(\delta_1-\delta^*_1)\wedge(\delta_2-\delta^*_3)$. Define
\[ s_0:= \inf\left\{S_{T_0}(\varphi):\varphi\in {\bf C}_{T_0},\varphi(0)\in \partial(\mathcal{A})_{\delta_2^*},
\rho_{T_0}\big(\varphi(\cdot),\Psi_{\cdot}(\varphi(0))\big)\geq \tilde{\delta} \right\}.  \]
Note that $s_0>0$ because of Hypothesis \ref{closeinfH}.

Proceeding as in the proof of Lemma \ref{klofR}, we assume that  $V(\partial (\mathcal{A})_{\delta_2},\partial (\mathcal{A})_{\delta_1})< s_0$. Then there exist a $T\in(0,+\infty)$ and $\psi\in {\bf C}_T, \psi(0)\in\partial(\mathcal{A})_{\delta_2},
\psi(T)\in\partial(\mathcal{A})_{\delta_1}$ such that
$S_T(\psi)<s_0$. Let $t^*:=\sup\{t\in[0,T]:\psi(t)\in \partial(\mathcal{A})_{\delta_2^*}\}$. Then by the continuity of $\psi$ and $\delta_2<\delta_2^*<\delta_1$, we get that $0<t^*<T$, $\psi(t^*)\in \partial(\mathcal{A})_{\delta_2^*}$ and $\psi(t)\notin \partial(\mathcal{A})_{\delta_2^*}, t\in (t^*, T].$ We extend the definition
of function $\psi$ up to the end of the interval $[0,t^*+T_0]$
as the solution of system {\rm(\ref{unpersys})} as long as $T < t^*+T_0$.  Thus $\psi(t^*+\cdot)|_{[0,T_0]}\in{\bf C}_{T_0}$. By the nonnegativity and additivity of $S$,
we have $S_{t^*,t^*+T_0}(\psi)\leq S_{t^*,T}\leq S_{T}(\psi)<s_0$.
From the definition of $s_0$, we get that
\begin{equation}\label{uses0agin}
\rho_{T_0}\big(\psi(t^*+\cdot),\Psi_{\cdot}(\psi(t^*))\big)< \tilde{\delta}.
\end{equation}

Case 1. $T_0\leq T-t^*$.  Combining (\ref{uses0agin}) with $\gamma^+\big(\Psi_{T_0}(\partial(\mathcal{A})_{\delta^*_2})\big)\subset (\mathcal{A})_{\delta^*_3}$, we know that $\psi(t^*+T_0)\in (\mathcal{A})_{\delta_2}$. Recall that $\psi(T)\in \partial(\mathcal{A})_{\delta_1}$ and use the continuity of $\psi$ again. Then there must be another time $\tilde{t}$ in $(t^*+T_0,T)$ such that $\psi(\tilde{t})\in \partial(\mathcal{A})_{\delta^*_2}$. This contradicts the definition of $t^*$.

Case 2. $T_0> T-t^*$.
For any $t\in[0,T_0]$, we obtain
\begin{align*}
  {\rm dist}\big(\psi(t^*+t),\mathcal{A}\big)
  \leq & \rho_{T_0}\big(\psi(t^*+t),\Psi_{t}(\psi(t^*))\big)+{\rm dist}\big(\Psi_{t}(\psi(t^*)),\mathcal{A}\big)\\
  < &\tilde{\delta}+\delta^*_1\\
  \leq &\delta_1.
\end{align*}
The second inequality above uses the fact that $\gamma^+\Big(\overline{(\mathcal{A})_{\delta^*_2}}\Big)\subset (\mathcal{A})_{\delta^*_1}$. This contradicts to the
condition $\psi(t^*+T-t^*)=\psi(T)\in \partial(\mathcal{A})_{\delta_1}$ while $T-t^*\in(0,T_0)$.
Thus we have $V(\partial(\mathcal{A})_{\delta_2},\partial(\mathcal{A})_{\delta_1})\geq s_0$.
\end{proof}

{\bf Proof of Theorem \ref{repeller}}.  Choose $\delta>0$ such that $(\mathcal{R})_{\delta}$ is a fundamental neighborhood of $\mathcal{R}$.
By Lemma \ref{klofR}, there exists $\delta_1,\ \delta_2\in(0,\delta),\ \delta_1>\delta_2$
and $s_0>0$ such that $V(\partial (\mathcal{R})_{\delta_1},\partial (\mathcal{R})_{\delta_2})\geq s_0$. Applying Proposition \ref{reppp} to $(\mathcal{R})_{\delta}$ and $\overline{(\mathcal{R})_{\delta_1}}\backslash(\mathcal{R})_{\delta_2}$, we have
$T_0=\inf\big\{u\geq 0:\Psi_{t}\big(\overline{(\mathcal{R})_{\delta_1}}\backslash(\mathcal{R})_{\delta_2}\big)
\subset((\mathcal{R})_{\delta})^c,t\geq u\big\}<+\infty$.
Let $F_0=\{\varphi\in {\bf C}_{T_0}:\varphi(0)\in \overline{(\mathcal{R})_{\delta_1}}\backslash (\mathcal{R})_{\delta_2}, \varphi(T_0)\in \overline{(\mathcal{R})_{\delta_1}}\}$
be a closed subset of ${\bf C}_{T_0}$. Then $F_0^0$ is bounded and $F_0$ does not contain any solution of system (\ref{unpersys}) by the definition of $T_0$.
By Hypothesis \ref{closeinfH},
we have $s_1=S_{T_0}(F_0)>0$.

 Let $\eta=\delta_2\wedge \frac{s_0\wedge s_1}{20L}\wedge \frac{s_0\wedge s_1}{20}>0$, where $L=L_{\overline{(\mathcal{R})_{\delta}}}$ is a constant as Lemma \ref{Vctin}. Then
by ${\bf(P_{\mathcal{R}})}$, there exist an open neighborhood $U$ of $\mathcal{R}$
and a positive constant $T^*$ such that for each
$x\in \mathcal{R}$
there are $\tilde{T}=\tilde{T}(x)\leq T^*,
\tilde{\varphi}^{x}\in {\bf C}_{\tilde{T}}$ satisfying $ \tilde{\varphi}^{x}(0)=x,
\tilde{\varphi}^{x}(\tilde{T})\in \partial(\mathcal{R})_{\delta_3}$ and $  S_{\tilde{T}}(\tilde{\varphi}^{x})<\eta$,
where $\delta_3=\eta\wedge\big({\rm dist}(\mathcal{R}, U^c)/2\big)>0$ so that $(\mathcal{R})_{\delta_3}\subset U$.
Here we have used the fact that the nonnegativity and additivity of $S$ and the continuity of $\tilde{\varphi}^{x}$.
By the compactness of $\overline{(\mathcal{R})_{\delta_1}}\backslash(\mathcal{R})_{\delta_3}$
in $ (\mathcal{R})_{\delta}\backslash\mathcal{R}$
and Proposition \ref{reppp}, $T_1=\inf\big\{u\geq 0:\Psi_{t}\big(\overline{(\mathcal{R})_{\delta_1}}\backslash(\mathcal{R})_{\delta_3}\big)
\subset((\mathcal{R})_{\delta})^c,t\geq u\big\}<+\infty$. Then obviously $T_0\leq T_1$.
For each $x\in
\overline{(\mathcal{R})_{\delta_2}}$,
define $\psi^x$ by
\[
\psi^x=\left\{
         \begin{array}{ll}
           \Psi_{\cdot}(x)|_{[0,T_1]}, & x\in \overline{(\mathcal{R})_{\delta_2}}\backslash(\mathcal{R})_{\delta_{3}}; \\
           {\rm LIF}_{x\hat{x}}\ast\tilde{\varphi}^{\hat{x}}\ast
\Psi_{\cdot}(\tilde{\varphi}^{\hat{x}}(\tilde{T}))|_{[0,T_1]}, & x\in (\mathcal{R})_{\delta_{3}},
         \end{array}
       \right.
\]
where $\hat{x}\in\mathcal{R}$ such that $|x-\hat{x}|={\rm dist}(x,\mathcal{R})<\delta_3$ if $x\in (\mathcal{R})_{\delta_3}$.
Let $T=T_1+\delta_3+T^*$.
Then we extend the domain of definition of $\psi^x$ up
to the end of the interval $[0,T]$ as the solution of {\rm(\ref{unpersys})},
 which does not increase the value of $S$. Thus for each $x\in\overline{(\mathcal{R})_{\delta_2}}$,
\begin{equation}\label{repup}
  S_{T}(\psi^x)\leq L\delta_3+\eta\leq 0.1(s_0\wedge s_1),
\end{equation}
and
\begin{equation}\label{reploc}
  \psi^x(T)\in((\mathcal{R})_{\delta})^c.
\end{equation}

We claim that there exists an $\varepsilon_*>0$ such that
\begin{equation}\label{repconout1}
  \int_{\big(\overline{(\mathcal{R})_{\delta_1}}\big)^c}
\mu^{\varepsilon}(dz)\mathbb{P}_z(X^{\varepsilon}_T\in (\mathcal{R})_{\delta_2})
\leq\exp\{-0.8s_0/\varepsilon^2\};
\end{equation}
\begin{equation}\label{repconout2}
  \int_{\overline{(\mathcal{R})_{\delta_1}}\backslash (\mathcal{R})_{\delta_2}}\mu^{\varepsilon}(dz)\mathbb{P}_z(X^{\varepsilon}_T\in(\mathcal{R})_{\delta_2})
\leq\exp\{-0.8(s_0\wedge s_1)/\varepsilon^2\};
\end{equation}
\begin{equation}\label{repconout3}
\int_{(\mathcal{R})_{\delta_2}}\mu^{\varepsilon}(dz)\mathbb{P}_z(X^{\varepsilon}_T\in(\mathcal{R})_{\delta_2})
\leq\mu^{\varepsilon}((\mathcal{R})_{\delta_2})-
\mu^{\varepsilon}((\mathcal{R})_{\delta_2})\exp\{-0.2(s_0\wedge s_1)/\varepsilon^2\},
\end{equation}
for any $\varepsilon\in(0,\varepsilon_*)$.

By the facts claimed above and the invariance of $\mu^{\varepsilon}$,
for any $\varepsilon\in(0,\varepsilon_*)$,
we have
\begin{align*}
\mu^{\varepsilon}((\mathcal{R})_{\delta_2})=&
\int_{\big(\overline{(\mathcal{R})_{\delta_1}}\big)^c}
\mu^{\varepsilon}(dz)\mathbb{P}_z(X^{\varepsilon}_T\in(\mathcal{R})_{\delta_2})
+\int_{\overline{(\mathcal{R})_{\delta_1}}\backslash (\mathcal{R})_{\delta_2}}\mu^{\varepsilon}(dz)\mathbb{P}_z(X^{\varepsilon}_T\in(\mathcal{R})_{\delta_2})
+\int_{(\mathcal{R})_{\delta_2}}\mu^{\varepsilon}(dz)\mathbb{P}_z(X^{\varepsilon}_T\in(\mathcal{R})_{\delta_2})\\
\leq& \exp\{-0.8s_0/\varepsilon^2\}+\exp\{-0.8(s_0\wedge s_1)/\varepsilon^2\}
+\mu^{\varepsilon}((\mathcal{R})_{\delta_2})-
\mu^{\varepsilon}((\mathcal{R})_{\delta_2})\exp\{-0.2(s_0\wedge s_1)/\varepsilon^2\}.
\end{align*}
Therefore, for any $\kappa\in(0,0.6(s_0\wedge s_1))$, shrinking $\varepsilon_*$ if necessary,
we have
\[   \mu^{\varepsilon}((\mathcal{R})_{\delta_2})\leq \exp\{-(0.8s_0-0.2(s_0\wedge s_1))/\varepsilon^2\}+\exp\{-0.6(s_0\wedge s_1)/\varepsilon^2\}\leq \exp\{-\kappa/\varepsilon^2\}
\:\ \textrm{ for all $\varepsilon\in(0,\varepsilon_*)$}.  \]

We now prove the claim. Define stopping times
\[ \eta_1=\eta_1^\varepsilon=\inf\{t\geq 0: X^{\varepsilon}_t\in\partial(\mathcal{R})_{\delta_1} \},\:
\eta_2=\inf\{t\geq \eta_1: X^{\varepsilon}_t\in\partial(\mathcal{R})_{\delta_2} \},\:
\tau=\inf\{t\geq 0: X^{\varepsilon}_t\in\partial(\mathcal{R})_{\delta_2} \}. \]
We shall prove the inequality
\begin{equation}\label{Esti1}
\mathbb{P}_z(X^{\varepsilon}_T \in (\mathcal{R})_{\delta_2} )\leq \sup_{y\in\partial(\mathcal{R})_{\delta_1}}\mathbb{P}_{y}(\tau\leq T),\forall z\in \big(\overline{(\mathcal{R})_{\delta_1}}\big)^c.
\end{equation}
In fact, for each $z\in \big(\overline{(\mathcal{R})_{\delta_1}}\big)^c$, by path continuity,
\[  \eta_2=\eta_1+\tau\circ\theta_{\eta_1}, \:\: a.s. \: \mathbb{P}_{z} .  \]
Thus, for each $z\in \big(\overline{(\mathcal{R})_{\delta_1}}\big)^c$, by path continuity and
the strong Markov property, we have
\begin{align*}
  \mathbb{P}_z(X^{\varepsilon}_T \in (\mathcal{R})_{\delta_2} )\leq &
\mathbb{P}_z(\eta_2 \leq T )\\
=&\mathbb{P}_z(\eta_1+\tau\circ\theta_{\eta_1} \leq T, \eta_1\leq T)\\
\leq&\mathbb{E}_{z}\big[1_{\{\eta_1\leq T\}} 1_{\{\tau\leq T\}}\circ\theta_{\eta_1} \big]\\
=&\mathbb{E}_{z}\big[ 1_{\{\eta_1\leq T\}}\mathbb{E}_{X^{\varepsilon}_{\eta_1}}[1_{\{\tau\leq T\}}] \big]\\
\leq& \sup_{y\in\partial(\mathcal{R})_{\delta_1}}\mathbb{P}_{y}(\tau\leq T).
\end{align*}
This proves (\ref{Esti1}).

Let $F_1=\{\varphi\in {\bf C}_T: \varphi(t) \in \partial(\mathcal{R})_{\delta_2}
\:\textrm{for some} \ t\in[0,T]\}$. Then $F_1$ is a closed subset of ${\bf C}_T$.
By the definition of $V(\partial(\mathcal{R})_{\delta_1},\partial(\mathcal{R})_{\delta_2})$
and Lemma \ref{klofR}, we have $\inf_{y\in \partial(\mathcal{R})_{\delta_1}}S_{T}^{y}(F_1)
\geq V(\partial(\mathcal{R})_{\delta_1},\partial(\mathcal{R})_{\delta_2})\geq s_0$.
By (\ref{ufldpupbb3}) of ${\bf(II_u')}$,
there exists $\varepsilon_1>0$ such that for any $\varepsilon\in(0,\varepsilon_1)$
and $y\in \partial(\mathcal{R})_{\delta_1}$, we get
\begin{align}\label{Esti2}
  \mathbb{P}_{y}(\tau\leq T)=\mathbb{P}_{y}(X^{\varepsilon}\in F_1)\leq
\exp\big\{-\big(\inf_{y\in \partial(\mathcal{R})_{\delta_1}}S_{T}^{y}(F_1)-0.1s_0\big)/\varepsilon^2\big\}
\leq \exp\{-0.9s_0/\varepsilon^2\}.
\end{align}
Therefore, by (\ref{Esti1}) and (\ref{Esti2}), for any $\varepsilon\in(0,\varepsilon_1)$,
\[\sup_{z\in \big(\overline{(\mathcal{R})_{\delta_1}}\big)^c}\mathbb{P}_z(X^{\varepsilon}_T \in (\mathcal{R})_{\delta_2} )\leq\sup_{y\in\partial(\mathcal{R})_{\delta_1}}\mathbb{P}_{y}(\tau\leq T)\leq \exp\{-0.9s_0/\varepsilon^2\}.
\]
Thus, the inequality (\ref{repconout1}) follows.

Let $F_2=\big\{\varphi\in {\bf C}_T: \varphi(0) \in \overline{(\mathcal{R})_{\delta_1}}\backslash (\mathcal{R})_{\delta_2}, \varphi(T) \in \overline{(\mathcal{R})_{\delta_2}}\big\}$.
Obviously, $F_2$ is a closed set in ${\bf C}_T$.
We claim that
\begin{equation}\label{Esti3}
S_{T}(F_2)\geq s_0\wedge s_1.
\end{equation}
Indeed, let $\tilde{\varphi}$ be any element in $F_2$ with $S_{T}(\tilde{\varphi})<s_1$. Then by the definition of $s_1$,
we have $\tilde{\varphi}(T_0)\notin \overline{(\mathcal{R})_{\delta_1}}$.
Thanks to the continuity of $\tilde{\varphi}$ and
$\tilde{\varphi}(T)\in \overline{(\mathcal{R})_{\delta_2}}$,
there exist $t_1,t_2\in(T_0,T],t_1<t_2$ such that
$\tilde{\varphi}(t_i)\in\partial(\mathcal{R})_{\delta_i},i=1,2$.
Therefore,  $S_{T}(\tilde{\varphi})\geq
S_{t_1t_2}(\tilde{\varphi})\geq V(\partial(\mathcal{R})_{\delta_1},\partial(\mathcal{R})_{\delta_2})
\geq s_0\geq  s_0\wedge s_1$. This proves the claim. By (\ref{ufldpupbb3}) of ${\bf(II_u')}$ and (\ref{Esti3}),
there exists $\varepsilon_2>0$ such that for any $\varepsilon\in(0,\varepsilon_2)$
and $z\in \overline{(\mathcal{R})_{\delta_1}}\backslash (\mathcal{R})_{\delta_2}$, we have
\[
\mathbb{P}_z(X^{\varepsilon}_T\in (\mathcal{R})_{\delta_2})\leq
\mathbb{P}_z(X^{\varepsilon}\in F_2)\leq
\exp\big\{-\big(\inf_{z\in \overline{(\mathcal{R})_{\delta_1}}\backslash
(\mathcal{R})_{\delta_2}}S_{T}^{z}(F_2)-0.1(s_0\wedge s_1)\big)/\varepsilon^2\big\}
\leq \exp\{-0.9(s_0\wedge s_1)/\varepsilon^2\}.
\]
Therefore the inequality (\ref{repconout2}) follows immediately.

Let $C=\{\varphi\in {\bf C}_T: \varphi(T)\notin (\mathcal{R})_{\delta_2}\},
G=\{\varphi\in {\bf C}_T: \rho_T(\varphi,\psi^{x})<\delta_*, x\in \overline{(\mathcal{R})_{\delta_2}}\}$,
where $\delta_*=\delta-\delta_2$. Then $G$ is an open set in ${\bf C}_T$
and $G\subset C$ by (\ref{reploc}).
Furthermore, $\sup_{z\in \overline{(\mathcal{R})_{\delta_2}}}S_{T}^{z}(G)\leq
\sup_{z\in \overline{(\mathcal{R})_{\delta_2}}}S_{T}(\psi^{z})\leq 0.1(s_0\wedge s_1)$ by (\ref{repup}).
By (\ref{ufldplbb3}) of ${\bf(I_u')}$,  there exists $\varepsilon_3>0$ such
that for any $\varepsilon\in(0,\varepsilon_3)$ and $z\in \overline{(\mathcal{R})_{\delta_2}}$,
we have
\begin{align*}
  \mathbb{P}_z(X^{\varepsilon}_T\notin (\mathcal{R})_{\delta_2})
=&\mathbb{P}_z(X^{\varepsilon}\in C)\\
\geq& \mathbb{P}_z(X^{\varepsilon}\in G)\\
\geq&\exp\big\{-\big(\sup_{z\in \overline{(\mathcal{R})_{\delta_2}}}S_{T}^{z}(G)+0.1(s_0\wedge s_1)\big)/\varepsilon^2\big\}\\
\geq&\exp\{-0.2(s_0\wedge s_1)/\varepsilon^2\}.
\end{align*}
 From this inequality, we obtain that for any $\varepsilon\in(0,\varepsilon_3)$ and $z\in \overline{(\mathcal{R})_{\delta_2}}$,
\begin{align*}
 \mathbb{P}_z(X^{\varepsilon}_T\in (\mathcal{R})_{\delta_2})
=&1- \mathbb{P}_z(X^{\varepsilon}_T\notin (\mathcal{R})_{\delta_2})\\
\leq &1-\exp\{-0.2(s_0\wedge s_1)/\varepsilon^2\}.
\end{align*}
So (\ref{repconout3}) follows easily from the above fact. This completes the proof under DZULDP. The proof under FWULDP is given in Appendix.

% The following property follows directly
%from the definition of topologically transitivity: for all open sets $U$ and $V$
%in $\mathbb{R}^r$ such that $U\cap \Lambda$ and $V\cap \Lambda$ are nonempty, there is $T>0$
%such that $\Psi_{T}(U)\cap V\neq\emptyset$.

In the following, we give three sufficient conditions to guarantee
the property ${\bf(P_{\mathcal{R}})}$ holds.
\begin{prp}\label{PR} A set $\mathcal{R}$ admits the property ${\bf(P_{\mathcal{R}})}$, if one of the following holds:

{\rm(i)} $\mathcal{R}$ is a compact set of Lebesgue measure zero;

{\rm(ii)} $\mathcal{R}$ is topologically transitive; and

{\rm(iii)} $\mathcal{R}$ is a nonempty compact set
satisfying $V(x,x_0)=0$ for all $x\in \mathcal{R}$ and some $x_0\in\mathcal{R}$,
and $V(x_0,z_0)=0$ for some $z_0\notin \mathcal{R}$.

\end{prp}
\begin{proof} (i) For every $\eta>0$, let $\delta=\frac{\eta}{2(L+1)}$ with $L=L_{\overline{(\mathcal{R})_{\eta}}}$ being a constant as Lemma \ref{Vctin}.
Since $\mathcal{R}$ is compact,
one can extract from the open cover $\bigcup_{x\in\mathcal{R}}B_{\delta}(x)$
of $\mathcal{R}$ a finite cover $\mathcal{R}$ by the sets $B_{\delta}(x_1),\cdots,B_{\delta}(x_N)$.
For each open set $B_{\delta}(x_k)$, since  $\mathcal{R}$ is a set of Lebesgue measure zero,
choose $y_k\in B_{\delta}(x_k)$ such that $y_k\notin \mathcal{R}$ for every $k=1,\cdots,N$.
Let $\delta_*:=\min_{1\leq k\leq N}{\rm dist}(y_k,\mathcal{R})$.
It is easy to see $\delta_*\in(0,\delta)$.
Take $U=(\mathcal{R})_{\delta_*}$ and $T^*=\eta<+\infty$.

For every $x\in\mathcal{R}$, since $\mathcal{R}\subset \bigcup_{k=1}^{N}B_{\delta}(x_k)$,
there exists $1\leq k\leq N$ such that $x\in B_{\delta}(x_k)$.
Let $T=T(x)=|x-x_k|+|x_k-y_k|$ and $\varphi^{x}={\rm LIF}_{xx_k}
\ast{\rm LIF}_{x_ky_k}$. Then $T<2\delta< \eta=T^*,\varphi^{x}\in {\bf C}_{T},\varphi^{x}(0)=x,\varphi^{x}(T)=y_k\notin U $,
and $S_T(\varphi^{x})\leq L|x-x_k|+L|x_k-y_k|<\eta$.
By definition, ${\bf(P_{\mathcal{R}})}$ holds.

 (ii) For every $\eta>0$, let $\delta=\frac{\eta}{4(L+1)}$ with $L=L_{\overline{(\mathcal{R})_{\eta}}}$ being a constant as Lemma \ref{Vctin}.
Note that $\mathcal{R}$ and $\partial(\mathcal{R})_{\delta}$ are  disjoint and compact.
Thus there exist $\hat{y}\in \mathcal{R}$ and $\hat{z}\in\partial(\mathcal{R})_{\delta}$ such
that $|\hat{y}-\hat{z}|={\rm dist}(\mathcal{R},\partial(\mathcal{R})_{\delta})=\delta$.
Since $\mathcal{R}$ is compact,
one can extract from the open cover $\bigcup_{x\in\mathcal{R}}B_{\delta}(x)$
of $\mathcal{R}$ a finite cover $\mathcal{R}$ by the sets $B_{\delta}(x_1),\cdots,B_{\delta}(x_N)$.
 Since both  sets $B_{\delta}(x_k)\cap\mathcal{R}$ and $B_{\delta}(\hat{y})\cap\mathcal{R}$ are nonempty, by the definition of
 topological transitivity, there exists $T_k>0$ such
that $\Psi_{T_k}(\hat{x}_k)\in B_{\delta}(\hat{y})$ for some $\hat{x}_k\in B_{\delta}(x_k)$,
where $k=1,\cdots,N$.

Let $U=(\mathcal{R})_{\delta}$ and $ T^*=\max_{1\leq k\leq N}T_k+4\delta<+\infty$.
For every $x\in\mathcal{R}$, there exists $1\leq k\leq N$ such that $x\in B_{\delta}(x_k)$.
Let $T=T(x)=|x-x_k|+|x_k-\hat{x}_k|+T_k+|\Psi_{T_k}(\hat{x}_k)-\hat{y}|+|\hat{y}-\hat{z}|$
and $\varphi^{x}={\rm LIF}_{xx_k}\ast{\rm LIF}_{x_k\hat{x}_k}\ast\Psi_{\cdot}(\hat{x}_k)|_{[0,T_k]}
\ast{\rm LIF}_{\Psi_{T_k}(\hat{x}_k)\hat{y}}\ast{\rm LIF}_{\hat{y}\hat{z}}$.
Then $T<4\delta+T_k\leq T^*,\varphi^{x}\in {\bf C}_{T},\varphi^{x}(0)=x,\varphi^{x}(T)=\hat{z}\notin U$,
and $S_T(\varphi^{x})\leq 4L\delta<\eta$.
By definition, ${\bf(P_{\mathcal{R}})}$ holds.

(iii) For every $\eta>0$, let $\delta=\frac{\eta}{4(L+1)}\wedge\big({\rm dist}(z_0,\mathcal{R})/2\big)$
and $U=(\mathcal{R})_{\delta}$, where $L=L_{\overline{(\mathcal{R})_{\eta}}}$
is a constant as Lemma \ref{Vctin}. Since $\mathcal{R}$ is compact,
one can extract from the open cover $\bigcup_{x\in\mathcal{R}}B_{\delta}(x)$
of $\mathcal{R}$ a finite cover $\bigcup_{k=1}^{N}B_{\delta}(x_k)$.
Because $V(x_k,x_0)=0$, there exist $T_k>0$ and $\varphi^k\in {\bf C}_{T_k}$ with $
\varphi^k(0)=x_k$ and $\varphi^k(T_k)=x_0$ such that
\[      S_{T_k}(\varphi^k)<\frac{\eta}{4} , \quad k=1,\cdots,N.     \]
Note that $V(x_0,z_0)=0$. Thus there exist $T_0>0$ and $\varphi^0\in {\bf C}_{T_0}$ with $
\varphi^0(0)=x_0$ and $\varphi^0(T_0)=z_0$ such that $S_{T_0}(\varphi^0)<\frac{\eta}{2}$.
Let $T^*=\delta+\max_{1\leq k\leq N}T_k+T_0<+\infty$. For every $x\in\mathcal{R}$, there exists $1\leq k\leq N$ such that $x\in B_{\delta}(x_k)$.
Let $T=T(x)=|x-x_k|+T_k+T_0$ and $ \varphi^{x}
={\rm LIF}_{xx_k}\ast \varphi^k\ast\varphi^0$.
Then $T< T^*,\varphi^{x}\in {\bf C}_T,\varphi^{x}(0)=x,\varphi^{x}(T)=z_0\in U^c$ and
\[  S_T(\varphi^{x})<L\delta+\frac{\eta}{4}+\frac{\eta}{2}<\eta. \]
This completes the proof.
\end{proof}

\begin{prp}\label{saddprp}
Suppose that there are points $x_1, x_2, \cdots, x_m$ in $\mathbb{R}^r$, an attractor $\mathcal{A}=:\mathcal{S}_{m+1}$ and compact equivalent classes $\mathcal{S}_1, \mathcal{S}_2, \cdots \mathcal{S}_m$  such that $\mathcal{A}$ and $\mathcal{S}_i\supset\alpha(x_i), i=1, 2,\cdots, m$ are pairwise disjoint, $\omega(x_i)\subset \mathcal{S}_{i+1}, i=1, 2,\cdots, m$.
Then $\bigcup_{i=1}^{m}\mathcal{S}_i\subset \Sigma:=\big\{y\in \mathbb{R}^r\backslash \mathcal{A}:V(y,\mathcal{A})=0\big\}$.
\end{prp}
\begin{proof}
We shall prove the conclusion by induction. Let $L=L_{\cup_{i=1}^{m+1}\overline{(\mathcal{S}_i)_{1}}}>1$ be a constant  given in Lemma \ref{Vctin}.

Assume that $m=1$. Then there exists $x_1\in\mathbb{R}^r$ such that $\alpha(x_1)\subset\mathcal{S}_1$ and $ \omega(x_1)\subset\mathcal{S}_{2}.$ Let $y\in \alpha(x_1)$ and $z\in \omega(x_1).$ For any $\eta\in(0,1)$, by the definition of limit point, there are $s_1>0$ and $t_1>0$ such that
$$T_1:=\big|\Psi_{-s_1}(x_1)-y\big|<\frac{\eta}{2L},\ T_2:=\big|\Psi_{t_1}(x_1)-z\big|<\frac{\eta}{2L}.$$
Define link function between $y$ and $z$ by
$$\psi^{yz}={\rm LIF}_{y\Psi_{-s_{1}}(x_{1})}\ast\Psi_{\cdot}\big(\Psi_{-s_1}(x_1)\big)|_{[0,t_1+s_1]}\ast{\rm LIF}_{\Psi_{t_1}(x_1)z}.$$
Then $S_T(\psi^{yz})<\eta$ with $T=T_1+T_2+s_1+t_1$. By the definition of $V(y,z)$, $V(y,z)=0$. For any $x\in \mathcal{S}_1$, we have
$$V(x,z)\le V(x,y)+V(y,z)=0$$
by the assumption that $\mathcal{S}_1$ is an equivalent class. Furthermore, $V(x,\mathcal{A})=0$, that is, $\mathcal{S}_1 \subset \Sigma$.

Suppose that the conclusion holds for $m-1$. Then $\bigcup_{i=2}^{m}\mathcal{S}_i\subset \Sigma$. Now choose $y\in \alpha(x_1)\subset\mathcal{S}_1$ and $z\in \omega(x_1)\subset\mathcal{S}_{2}.$ Proceeding as in the case of $m=1$, we get that $V(y,z)=0$. Thus, $V(y,\mathcal{A})\le V(y,z)+V(z,\mathcal{A})=0$ by $V(y,z)=0$ and the induction assumption. Similarly, $V(x,\mathcal{A})=0$ for any $x\in \mathcal{S}_1$. This proves $\mathcal{S}_1 \subset \Sigma$. Together with the induction assumption, we conclude that $\bigcup_{i=1}^{m}\mathcal{S}_i\subset \Sigma$.
\end{proof}

From Proposition \ref{saddprp}, Theorem \ref{saddle} is a corollary of the following theorem.
\begin{theorem}\label{mainatt}
Let the diffusion matrix $a$ be nonsingular on $\mathbb{R}^r$ and system {\rm(\ref{itodff})} satisfy Hypothesis \ref{closeinfH} and possess {\rm DZULDP} or {\rm FWULDP}.  Suppose that
 $\mathcal{A}$ is an attractor and $\mu^{\varepsilon}$ is any stationary distributions of diffusion process
$(X^{\varepsilon},\mathbb{P}_x)$. If $y\in \mathbb{R}^r\backslash \mathcal{A}$ such that $V(y,\mathcal{A})=0$,
then there exist $\delta=\delta(y)>0$, $\kappa>0$ and $\varepsilon_*>0$
such that for any $\varepsilon\in(0,\varepsilon_*)$ we have
\begin{equation}\label{Expdec}
 \mu^{\varepsilon}\big(B_{\delta}(y)\big)\leq \exp\{-\kappa/ \varepsilon^2\}  .   \end{equation}
As a result, if, moreover, $\mu^{\varepsilon}_j\xlongrightarrow{w} \mu$ as $\varepsilon_j\rightarrow 0$,
then $\mu(B_{\delta}(y))=0$.
In particular, $\mu(\Sigma)=0$ with $\Sigma:=\{y\in\mathbb{R}^r\backslash \mathcal{A}:V(y,\mathcal{A})=0\}$. Besides, if $K\subset \Sigma$ is a nonempty compact subset, then there is an open set $U\supset K$, $\tilde{\kappa}>0$ and $\tilde{\varepsilon_*}>0$
such that for any $\varepsilon\in(0,\tilde{\varepsilon_*})$ we have
\begin{equation}\label{Expdec1}
 \mu^{\varepsilon}\big(U\big)\leq \exp\{-\tilde{\kappa}/ \varepsilon^2\}  .   \end{equation}
\end{theorem}
\begin{proof}
Let $\delta'\in(0,{\rm dist}\left(y,\mathcal{A}\right)/2)$ such that
$(\mathcal{A})_{\delta'}$ is a fundamental neighborhood of $\mathcal{A}$.
By Lemma \ref{klofA}, there exist $\delta_1>\delta_2 >\delta_3,\ \delta_1,\delta_2, \delta_3 \in(0,\delta')$
and $s_0>0$ such that $V\left(\partial (\mathcal{A})_{\delta_2},\partial (\mathcal{A})_{\delta_1}\right)\geq s_0$.
Since $\mathcal{A}$ is an attractor,  we obtain
$T_0:=\inf\big\{u\geq 0:\Psi_{t}\big(\overline{(\mathcal{A})_{\delta_1}}\big)\subset(\mathcal{A})_{\delta_3}, t\geq u\big\}<+\infty$. The set $F_0=\{\varphi\in {\bf C}_{T_0}:\varphi(0)\in \overline{(\mathcal{A})_{\delta_1}},\varphi(T_0)\in ((\mathcal{A})_{\delta_2})^c\}$ is a closed subset of ${\bf C}_{T_0}$, $F_0^0$ is bounded and $F_0$ does not contain any solution of system (\ref{unpersys}).
Thus, by Hypothesis \ref{closeinfH}, $s_1:=S_{T_0}(F_0)>0$.
Let $\delta=\delta_2\wedge \frac{s_0\wedge s_1}{10L}$,
where $L=L_{\bar{B}_{\delta'}(y)\cup \overline{(\mathcal{A})_{\delta'}}}$ is a constant as Lemma \ref{Vctin}.
Since $V(y,\mathcal{A})=0$, there exist $T_1>0$ and $ \tilde{\varphi}\in {\bf C}_{T_1}$ with $
\tilde{\varphi}(0)=y$ and $\tilde{\varphi}(T_1)\in \mathcal{A}$ such that $S_{T_1}(\tilde{\varphi})<0.2(s_0\wedge s_1)$.

For each $x\in \bar{B}_{\delta}(y)$, define $\psi^x$ as follows:
\[ \psi^x={\rm LIF}_{xy}\ast \tilde{\varphi}.\]
Note that the domain of $\psi^x$ is the subinterval of $[0,\delta+T_1]$.
Let $T=T_0\vee(\delta+T_1)$.
Then we extend the definition domain of $\psi^x$ up
to the end of the interval $[0,T]$ as the solution of (\ref{unpersys}),
which does not increase the value of $S$. Thus
\begin{equation}\label{attractup}
 S_{T}(\psi^x)\leq L\delta+0.2(s_0\wedge s_1)\leq 0.3(s_0\wedge s_1)
\end{equation}
for every $x\in\bar{B}_{\delta}(y)$.
By the invariance of $\mathcal{A}$, we also have
\begin{equation}\label{attraloc}
\psi^x(T)\in \mathcal{A}\ {\rm for\ every}\ x\in\bar{B}_{\delta}(y).
\end{equation}

We claim that there exists $\varepsilon_*>0$ such that
\begin{equation}\label{saddleep10}
\int_{B_{\delta}(y)}\mu^{\varepsilon}(dz)\mathbb{P}_z\big(X^\varepsilon_T\in (\mathcal{A})_{\delta_1}\big)\geq\mu^{\varepsilon}\big(B_{\delta}(y)\big)\exp\{-0.4(s_0\wedge s_1)/\varepsilon^2\};\ {\rm and}
\end{equation}
\begin{equation}\label{saddleep20}
  \int_{(\mathcal{A})_{\delta_1}}\mu^{\varepsilon}(dz)\mathbb{P}_z\big(X^\varepsilon_T\in (\mathcal{A})_{\delta_1}\big)
\geq \mu^{\varepsilon}((\mathcal{A})_{\delta_1})-\mu^{\varepsilon}((\mathcal{A})_{\delta_1})
\exp\{-0.9(s_0\wedge s_1)/\varepsilon^2\}
\end{equation}
for each $\varepsilon\in(0,\varepsilon_*)$.

By (\ref{saddleep10}), (\ref{saddleep20}) and the invariance of $\mu^{\varepsilon}$,
for any $\varepsilon\in(0,\varepsilon_*)$,
we have
\begin{align*}
  \mu^{\varepsilon}((\mathcal{A})_{\delta_1})=&\int_{((\mathcal{A})_{\delta_1})^c}\mu^{\varepsilon}(dz)\mathbb{P}_z\big(X^\varepsilon_T\in (\mathcal{A})_{\delta_1}\big)+\int_{(\mathcal{A})_{\delta_1}}\mu^{\varepsilon}(dz)\mathbb{P}_z\big(X^\varepsilon_T\in (\mathcal{A})_{\delta_1}\big)\\
\geq&
\int_{B_{\delta}(y)}\mu^{\varepsilon}(dz)\mathbb{P}_z\big(X^\varepsilon_T\in (\mathcal{A})_{\delta_1}\big)
+\int_{(\mathcal{A})_{\delta_1}}\mu^{\varepsilon}(dz)\mathbb{P}_z\big(X^\varepsilon_T\in (\mathcal{A})_{\delta_1}\big)\\
\geq&\mu^{\varepsilon}\big(B_{\delta}(y)\big)\exp\{-0.4(s_0\wedge s_1)/\varepsilon^2\}+\mu^{\varepsilon}((\mathcal{A})_{\delta_1})-\mu^{\varepsilon}((\mathcal{A})_{\delta_1})
\exp\{-0.9(s_0\wedge s_1)/\varepsilon^2\}.
\end{align*}
Therefore,
\[ \mu^{\varepsilon}\big(B_{\delta}(y)\big)\leq
\mu^{\varepsilon}(\mathcal{A}_{\delta_1})\exp\{-0.5(s_0\wedge s_1)/\varepsilon^2\}
\leq \exp\{-\kappa/\varepsilon^2\},   \]
for each $\varepsilon\in(0,\varepsilon_*)$, where $\kappa=0.5(s_0\wedge s_1)$.
(\ref{Expdec}) has been proved.

To prove the claim, let
$G=\{\varphi\in {\bf C}_{T}:\varphi(T)\in (\mathcal{A})_{\delta_1}\}$
and $F=\{\varphi\in {\bf C}_{T}:\varphi(T)\notin (\mathcal{A})_{\delta_1}\}$.
Obviously, $G$ is open and $F$ is closed in ${\bf C}_{T}$.

For any $z\in \bar{B}_{\delta}(y)$, by (\ref{attraloc}),
we have $\psi^{z}(T)\in \mathcal{A}\subset (\mathcal{A})_{\delta_1}$.
Thus $\psi^{z}\in G$. This implies that
\[   \sup_{z\in \bar{B}_{\delta}(y)}S_{T}^{z}(G)\leq \sup_{z\in \bar{B}_{\delta}(y)} S_{T}(\psi^z)\leq 0.3(s_0\wedge s_1)\ {\rm by}\ (\ref{attractup}).   \]
By (\ref{ufldplbb3}) of ${\bf(I_u')}$,  there exists $\varepsilon_1>0$ such that for any $\varepsilon\in(0,\varepsilon_1),z\in\bar{B}_{\delta}(y)$ we have
\begin{equation}\label{saddleep110}
  \mathbb{P}_{z}\big(X^{\varepsilon}_T\in (\mathcal{A})_{\delta_1} \big)
=\mathbb{P}_{z}\big(X^{\varepsilon}\in G\big)\geq \exp\big\{-\big(0.3(s_0\wedge s_1)+0.1(s_0\wedge s_1)\big)/\varepsilon^2\big\}
=\exp\{-0.4(s_0\wedge s_1)/\varepsilon^2\}.
\end{equation}
Hence (\ref{saddleep10}) follows from (\ref{saddleep110}).

Next we show that $\inf_{z\in\overline{(\mathcal{A})_{\delta_1}}}S^{z}_T(F)\geq (s_0\wedge s_1)$. In fact, let $\tilde{\varphi}$ be any element in $F$ satisfying $\tilde{z}:=\tilde{\varphi}(0)\in \overline{(\mathcal{A})_{\delta_1}}$.
If $S_{T}^{\tilde{z}}(\tilde{\varphi})<s_1$,
then $S_{T_0}^{\tilde{z}}(\tilde{\varphi})<s_1$. So $\tilde{\varphi}(T_0)\in (\mathcal{A})_{\delta_2}$ from the definition of $s_1$.
Note that $\tilde{\varphi}(T)\notin (\mathcal{A})_{\delta_1}$
and $\overline{(\mathcal{A})_{\delta_2}}\subset (\mathcal{A})_{\delta_1}$.
Then by the continuity of $\tilde{\varphi}$, there exist $T_1,T_2\in(T_0,T],T_1>T_2$
such that $\tilde{\varphi}(T_i)\in \partial (\mathcal{A})_{\delta_i},i=1,2$.
Thus $S_{T}^{\tilde{z}}(\tilde{\varphi})\geq
S_{T_2T_1}^{\tilde{z}}(\tilde{\varphi})\geq
V(\partial (\mathcal{A})_{\delta_2},\partial (\mathcal{A})_{\delta_1})\geq s_0$.
Therefore $\inf_{z\in \overline{(\mathcal{A})_{\delta_1}}}S^{z}_T(F)\geq (s_0\wedge s_1)$.
Finally, by (\ref{ufldpupbb3}) of ${\bf(II_u')}$,
there exists $\varepsilon_2>0$ such that for any
$\varepsilon\in(0,\varepsilon_2),z\in \overline{(\mathcal{A})_{\delta_1}}$ we have
\[\mathbb{P}_z\big(X^\varepsilon_T\notin (\mathcal{A})_{\delta_1}\big)
=\mathbb{P}_z\big(X^\varepsilon\in F\big)\leq
\exp\big\{-\big(\inf_{z\in\overline{(\mathcal{A})_{\delta_1}}}S^{z}_T(F)-0.1(s_0\wedge s_1)\big)/\varepsilon^2\big\}
\leq\exp\{-0.9(s_0\wedge s_1)/\varepsilon^2\} . \]
This implies
\begin{align*}
   \int_{(\mathcal{A})_{\delta_1}}\mu^{\varepsilon}(dz)\mathbb{P}_z\big(X^\varepsilon_T\in (\mathcal{A})_{\delta_1}\big)=&
 \int_{(\mathcal{A})_{\delta_1}}\mu^{\varepsilon}(dz)\big(1-\mathbb{P}_z\big(X^\varepsilon_T\notin (\mathcal{A})_{\delta_1}\big)\big)\\
=&\mu^{\varepsilon}((\mathcal{A})_{\delta_1})
-\int_{(\mathcal{A})_{\delta_1}}\mu^{\varepsilon}(dz)\mathbb{P}_z\big(X^\varepsilon_T\notin (\mathcal{A})_{\delta_1}\big)\\
\geq& \mu^{\varepsilon}((\mathcal{A})_{\delta_1}) -\mu^{\varepsilon}((\mathcal{A})_{\delta_1})
\exp\{-0.9(s_0\wedge s_1)/\varepsilon^2\},
\end{align*}
for any $\varepsilon\in(0,\varepsilon_2)$, i.e., the inequality (\ref{saddleep20}) follows.

Furthermore, if $\mu^{\varepsilon}_j\xlongrightarrow{w} \mu$ as $\varepsilon_j\rightarrow 0$,
then $\mu\big(B_{\delta}(y)\big)\leq\liminf_{\varepsilon_j\rightarrow 0}\mu^{\varepsilon_j}\big(B_{\delta}(y)\big)=0$ from (\ref{Expdec}).
In particular, for every $y\in\Sigma$, there exists $\delta(y)>0$
such that $\mu(B_{\delta}(y))=0$.
Note that $\{B(y,\delta(y))\}_{y\in\Sigma}$ is an open cover of $\Sigma$,
by the Lindel\"{o}f theorem, there exists a countable subcover $\{B(y_k,\delta(y_k))\}_{k\in\mathbb{N}}$ of $\Sigma$.
Therefore, by the subadditivity of $\mu$, $\mu(\Sigma)\leq \sum_{k=1}^{\infty}\mu\big(B(y_k,\delta(y_k))\big)=0$.

Let $K\subset \Sigma$ be compact. Then from the open cover $\{B(y,\delta(y))\}_{y\in K}$ of $K$ we can extract a finite open cover $\bigcup_{i=1}^{N}B\big(y_i,\delta(y_i)\big)\supset K\supset\{y_1, y_2,\cdots,y_N\}$. Applying (\ref{Expdec}) to $y_1, y_2,\cdots,y_N$, we obtain that there exist $\kappa_i>0$ and $\varepsilon_*^i>0$
such that for any $\varepsilon\in(0,\varepsilon_*^i)$ we have
\begin{equation}\label{Expdeck}
 \mu^{\varepsilon}\big(B_{\delta(y_i)}(y_i)\big)\leq \exp\{-\kappa_i/ \varepsilon^2\},\ i=1,2,\cdots,N.
 \end{equation}
 Denote by $U$ the open cover $\bigcup_{i=1}^{N}B\big(y_i,\delta(y_i)\big)$. Let $\tilde{\kappa}\in \big(0, (\kappa_1\wedge \kappa_2\wedge\cdots \wedge \kappa_N)\big)$. Then there exists $\tilde{\varepsilon_*}\in \big(0, (\varepsilon_*^1\wedge \varepsilon_*^2\wedge\cdots \wedge \varepsilon_*^N)\big)$ such that (\ref{Expdec1}) holds for any $\varepsilon\in(0,\tilde{\varepsilon_*})$. The proof under DZULDP is complete. The proof under FWULDP is given in Appendix.
\end{proof}

\section{Applications}
In this section, we first apply our main results to the Morse--Smale system, the Axiom A system, those systems for the Poincar\'{e}--Bendixson property to hold and gradient systems or gradient-like systems and then present a number of nontrivial examples whose limit measures are precisely obtained, hence their concentrations are accurately depicted, theses examples contain monostable systems, multistable systems and those having infinite equivalent classes.

It is well known that the Poincar\'{e}--Bendixson theorem holds for planar system {\rm(\ref{unpersys})}, that is, any limit set for a planar system containing no equilibria is a periodic orbit, see \cite {Hartman}. The Poincar\'{e}--Bendixson theorem still holds for some higher dimensional systems, for example, three dimensional competitive and cooperative systems
\cite{Hirsch, Smith}; monotone systems with cones of rank two \cite{WangYi, Sanchez}; monotone cyclic feedback systems \cite{SmithJ, SellJ} and many references therein.

A point $x$ is {\it nonwandering} for $\Psi$ if for every neighborhood $U$ of $x$, and $T>0$, there is $t>T$ such that $\Psi_t(U)\cap U\neq\emptyset$. The set of nonwandering points of $\Psi$ is denoted by $\Omega(\Psi)$ (or simply $\Omega$
if there is no risk of confusion). It is easy to see that $\Omega$ contains all limit points for $\Psi$, that is, $\Omega \supset L$ where $L=L(\Psi)$ denotes the set of all limit points for $\Psi$.

The deterministic system {\rm(\ref{unpersys})} is  called {\it Morse--Smale} \cite{Palis,Palis1, Palis2, Benaim} if

(i) {\rm(\ref{unpersys})} has a finite number of critical elements (equilibria and periodic orbits) all of which are hyperbolic;

(ii) if $\sigma_1$ and $\sigma_2$ are critical elements of {\rm(\ref{unpersys})} then $W^s(\sigma_1)$ is transversal to
$W^u(\sigma_2)$;

(iii) {\rm(\ref{unpersys})} admits a global attractor $\mathcal{A}$; and

(iv) the nonwandering set is equal to the union of the critical elements of $\Psi$.

Note that  from Theorem 3.2 of Wilson \cite{Wilson} it follows that any attractor $\mathcal{A}$ has a $C^{\infty}$ Liapunov function $V: \mathbb{R}^r\rightarrow \mathbb{R}_+$ satisfying $\lim_{|x|\rightarrow +\infty}V(x)=+\infty$ and  $\langle b(x),\nabla V(x)\rangle <0$ for $x\in \mathcal{A}^c$. This means that $b$ is transversal to
$V^{-1}(R)$ for $R$ sufficiently large.

Let $B(\Psi)=\overline{\{x\in \mathbb{R}^r: x\in \omega(x)\}}$ be the Birkhoff center of $\Psi$. Recall from \cite{Chen2} that ${\rm supp}(\mu)\subset B(\Psi)$ for any limit measure $\mu$ of $\{\mu^{\varepsilon}\}$ in the weak*-topology.

For a given compact invariant set $K\subset \mathbb{R}^r$, we  define
 $$W^s(K)=\big\{x\in \mathbb{R}^r\big|\lim_{t\rightarrow +\infty}{\rm dist}\big(\Psi_t(x),K\big)=0\big\}\ {\rm and}$$
 $$W^u(K)=\big\{x\in \mathbb{R}^r\big|\lim_{t\rightarrow -\infty}{\rm dist}\big(\Psi_t(x),K\big)=0\big\}$$
as the stable set and the unstable set of $K$, respectively.
\begin{lemma}\label{AttractorCri}
Let $K\subset \mathbb{R}^r$ be a compact invariant set for $\Psi$ and $N$ a neighborhood of $K$ satisfying that $\overline{N}$ is compact and $\overline{N}\cap L =K$. Then $K$ is an attractor {\rm(}repeller{\rm)} if and only if $W^u(K)\backslash K=\emptyset$ {\rm(}$W^s(K)\backslash K=\emptyset${\rm)}. The same result still holds when $L$ is replaced by $\Omega$.
\end{lemma}

 \begin{proof}
We only prove the sufficiency for attractor. Suppose the contrary. Then $K$ is not an attractor. Let $A$ be the maximal invariant set of $\Psi$ in $\overline{N}$. Then $A\supset K$. We claim that $A=K$. Otherwise, take $x\in A\backslash K$. Then by the invariance of $A$, $\gamma(x)\subset A$ and furthermore $\alpha(x)\subset \overline{N}\cap L=K$.  This implies that $x\in W^u(K)\backslash K$, a contradiction.   By \cite[Lemma 5.4, p.72]{Benaim1}, there exists a $p\in \partial N$ such that $\gamma^{-}(p)\subset N$. Thus $\alpha(p)\subset \overline{N}\cap L =K$, hence  $p\in W^u(K)\backslash K$, a contradiction to $W^u(K)\backslash  K=\emptyset$. The proof is still valid when $L$ is replaced by $\Omega$.

 \end{proof}
\begin{cor} \label{cor}
Let $A$ be nonsingular on $\mathbb{R}^r$, system {\rm(\ref{itodff})} possess {\rm DZULDP} or {\rm FWULDP} and system {\rm(\ref{unpersys})} have a finite number of equilibria and periodic orbits.
Suppose that the system {\rm(\ref{unpersys})} admits the Poincar\'{e}--Bendixson property without any cycle.
Let $p_1,\cdots,p_m$ and $\gamma_1,\cdots,\gamma_l$ denote all the asymptotically stable equilibria
and all the asymptotically stable periodic orbits of {\rm(\ref{unpersys})}, respectively.
Then every limit measure $\mu$ of $\{\mu^{\varepsilon}\}$ must be  a convex combination of the measures
\[    \delta_{p_1}(\cdot),\cdots,\delta_{p_m}(\cdot), \frac{1}{T_1}\int_{0}^{T_1}\delta_{\Psi_{t}(x_1)}(\cdot)dt,
\cdots, \frac{1}{T_l}\int_{0}^{T_l}\delta_{\Psi_{t}(x_l)}(\cdot)dt,   \]
where $T_k$ is the period of $\gamma_k$ and $x_k$ is any point in $\gamma_k, k=1,\cdots,l$.
\end{cor}
\begin{proof}
Since $\Psi$ admits the Poincar\'{e}--Bendixson property and has no cycle. $L$ is the union of all equilibria and closed orbits. Let $\sigma\subset L$ be an equilibrium or a closed orbit. Then ${\rm dist}(\sigma,L\backslash \sigma)>0$ by the assumptions. Take $\delta\in (0,0.5{\rm dist}(\sigma,L\backslash \sigma))$ and $N=(\sigma)_{\delta}$. Then $\overline{N}\cap L=\sigma$. Applying Lemma \ref{AttractorCri}, we obtain that $\sigma$ is an attractor (a repeller) if and only if $W^u(\sigma)\backslash \sigma=\emptyset$ ($W^s(\sigma)\backslash\sigma=\emptyset$). By Theorem \ref{repeller}, $\mu(\sigma)=0$ if $W^s(\sigma)\backslash\sigma=\emptyset$.

Let $\sigma_1\subset L$ such that both $W^s(\sigma_1)\backslash\sigma_1$ and $W^u(\sigma_1)\backslash\sigma_1$ are nonempty. Take $x_1\in W^u(\sigma_1)\backslash \sigma_1 $. Then $\omega(x_1)=:\sigma_2\subset L$ with $\sigma_2\neq \sigma_1$ by no-cycles assumption. If $W^u(\sigma_2)\backslash\sigma_2=\emptyset$, then $\sigma_2$ is an attractor. By Theorem \ref{saddle}, $\mu(\sigma_1)=0$. Otherwise, there exists a point $x_2\in W^u(\sigma_2)\backslash \sigma_2 $. Continuing this procedure, it must end in a finite step, say $m-$th step, because of the finiteness of equilibria and closed orbits. This means that we have $\{\sigma_1,\sigma_2,\cdots, \sigma_m,\sigma_{m+1} \}\subset L$ and $x_i\in \big(W^u(\sigma_i)\backslash \sigma_i\big)\bigcap \big(W^s(\sigma_{i+1})\backslash \sigma_{i+1}\big)$ for $i=1,2,\cdots, m$ such that $\alpha(x_i)=\sigma_i$ and $\omega(x_i)=\sigma_{i+1}$ for $i=1,2,\cdots, m$ and $\sigma_{m+1}$ is an attractor. Applying Theorem \ref{saddle}, we get that $\mu(\sigma_i)=0$ for $i=1,2,\cdots, m$. Since the Birkhoff center of {\rm(\ref{unpersys})} is $B(\Psi)=L$ under our assumptions, the conclusion follows immediately from \cite[Theorem 2.1]{Chen2}.
\end{proof}

Applying Corollary \ref {cor} to Examples 4.2, 4.4 and 4.6 and two freedom degree system in Remark 4.10 in \cite{Chen2020} (please refer to it), we get rid of any concentration on saddles and repellers and proves the conjecture proposed there. Corollary \ref {cor} can be applied to not only Morse--Smale systems but also structural unstable systems, see Proposition \ref{Van der Pole} below.

{\rm(\ref{unpersys})} is called a {\it gradient--like} system if it admits a Liapounov function which strictly decreases along nonconstant forward orbits of {\rm(\ref{unpersys})}. The following result improves Theorem 2.1 of {\rm\cite{Huang2016}}.
\begin{cor}\label{gradient}
 Suppose that {\rm(\ref{unpersys})} is a gradient system or a gradient--like system which has a finite number of equilibria. Then any limit measure of the corresponding system {\rm(\ref{itodff})} is concentrated on Liapunov stable equilibria.
 \end{cor}
 \begin{proof}
 Let $b(x)=-\nabla F(x)$ for a $C^2$ function $F$ with finite critical points. We claim that {\rm(\ref{unpersys})} has no cycle of connecting orbits. Suppose the contrary.  Then there exist critical points $P_1, P_2, \cdots, P_n$ and points $x_1, x_2, \cdots, x_{n}$ such that $\alpha(x_i)=P_i$ and $\omega(x_i)=P_{i+1}$ for $1, 2, \cdots, n$ with $P_{n+1}=P_1$. Since $\frac{d}{dt}F\big(\Psi_t(x_i)\big)=-|\nabla F\big(\Psi_t(x_i)\big)|^2<0$ for each $i$, we get that $F(P_1)<F(P_2)<\cdots <F(P_n)<F(P_1)$,  a contradiction. Thus we conclude that any limit measure of the corresponding system {\rm(\ref{itodff})} is concentrated away from all maximal points and saddle points. The same result for a gradient--like system can be proved in the same manner.
 \end{proof}

Let $F(x_1,x_2,x_3)=x_1^2(x_1-1)^2+x_2^2(x_2-1)^2+x_3^2$ or $F(x_1,x_2,\cdots,x_n)=\sum_{i=1}^nx_i^2(x_i-1)^2$. Then From Corollary \ref{gradient} it follows that limit measures  are supported on the minimal points with global minimum potential, which coincides with the well--known result for gradient systems perturbed by additive noise, see \cite{Hwang1980}.

Assume that {\rm(\ref{unpersys})} is dissipative. Then as done in Bena\"{\i}m \cite {Benaim1},  we can always suppose (by multiplying $b$ by a
smooth positive function which goes to zero as $|x|\rightarrow +\infty$) that the flow induced by $b$
is defined on the $r$-sphere $S^r= \mathbb{R}^r\cup\{\infty\}$, where the point at infinity is a source.

{\rm(\ref{unpersys})} is called {\it Axiom A} \cite {Smale}, if its nonwandering set $\Omega$ is hyperbolic and its critical elements are dense in $\Omega$. According to the spectral decomposition theorem \cite {Smale}, there is a unique way of writing
$\Omega$ as the finite union of disjoint, closed, invariant indecomposable subsets on each of which $\Psi$ is topologically transitive:
$$\Omega=\Omega_1\cup\Omega_2\cup\cdots \cup\Omega_k.$$
Such $\Omega_i\ (i=1,\cdots,k)$ are called the {\it  basic sets} of $\Psi$.

Now we define the relation $\succ$ on the basic sets:
$$\Omega_i\succ \Omega_j\Longleftrightarrow \big(W^u(\Omega_i)\backslash \Omega_i\big)\cap \big(W^s(\Omega_j)\backslash \Omega_j\big)\neq \emptyset.$$
$\{\Omega_{i_1},\Omega_{i_2},\cdots, \Omega_{i_m}\}$ with $1\le m\le k$ is called an {\it $m-$cycle} if
$$
\Omega_{i_1}\succ\Omega_{i_2}\succ\cdots\succ \Omega_{i_m}\succ \Omega_{i_1}.
$$
\begin{cor}\label{AxiomA}
 Suppose that {\rm(\ref{unpersys})} is an Axiom A system.  Then any limit measure of the corresponding system {\rm(\ref{itodff})} is concentrated on Liapunov stable basic sets or cycles. In particular, if {\rm(\ref{unpersys})} is $\Omega-$stable, or structurally  stable, or a separated star, then any limit measure  is concentrated on Liapunov stable basic sets.
 \end{cor}
\begin{proof}
Let $\mu^{\varepsilon}_j\xlongrightarrow{w} \mu$ as $\varepsilon_j\rightarrow 0$. For a given basic set $\Omega_i$, choose $\delta\in (0,0.5{\rm dist}(\Omega_i,\Omega\backslash \Omega_i))$ and $N=(\Omega_i)_{\delta}$. Then $\overline{N}\cap \Omega=\Omega_i$. Applying Lemma \ref{AttractorCri}, we obtain that $\Omega_i$ is an attractor (a repeller) if $W^u(\Omega_i)\backslash \Omega_i=\emptyset$ ($W^s(\Omega_i)\backslash \Omega_i=\emptyset$). By Theorem \ref{repeller}, $\mu(\Omega_i)=0$ if $W^s(\Omega_i)\backslash\Omega_i=\emptyset$.

Suppose that $\Omega_{i_1}$ is Liapunov unstable. Then $W^u(\Omega_{i_1})\backslash \Omega_{i_1}\neq\emptyset$.  Choose $x_1\in W^u(\Omega_{i_1})\backslash \Omega_{i_1}$. Thanks to the invariance and non-intersection of the basic sets and the definition of $W^u(\Omega_{i_1})$,  $x_1\notin \Omega$. Since $\omega(x_1)\subset \Omega$ and is connected,  there exists an index $i_2\in I:=\{1,2,\cdots,k\}$ such that $\omega(x_1)\subset \Omega_{i_2}$, that is, $x_1\in W^s(\Omega_{i_2})\backslash \Omega_{i_2}$. By definition, $\Omega_{i_1}\succ \Omega_{i_2}$. If $i_2=i_1$, then $\Omega_{i_1}$ forms a $1$-cycle and the conclusion is true. Otherwise, $i_2\neq i_1$. Suppose that $W^u(\Omega_{i_2})\backslash \Omega_{i_2}=\emptyset$. Then $\Omega_{i_2}$ is Liapunov stable. It follows from the topological transitivity of basic sets and Theorem \ref{saddle} that $\mu(\Omega_{i_1})=0$. Let $W^u(\Omega_{i_2})\backslash \Omega_{i_2}\neq \emptyset$. Then there are an $x_2\in W^u(\Omega_{i_2})\backslash \Omega_{i_2}$ and an index $i_3\in I$ such that $x_2\in W^s(\Omega_{i_3})\backslash \Omega_{i_3}$. So far we have $\Omega_{i_1}\succ \Omega_{i_2}\succ \Omega_{i_3}$. If $i_3=i_1$, then $\{\Omega_{i_1}, \Omega_{i_2}\}$ forms a $2$-cycle and the conclusion holds. Otherwise, $i_3\neq i_1$. Suppose that $W^u(\Omega_{i_3})\backslash \Omega_{i_3}=\emptyset$. Then $\Omega_{i_3}$ is Liapunov stable. It follows from the topological transitivity of basic sets and Theorem \ref{saddle} that $\mu(\Omega_{i_j})=0$ for $j=1,2$. Let $W^u(\Omega_{i_3})\backslash \Omega_{i_3}\neq \emptyset$. Then there are an $x_3\in W^u(\Omega_{i_3})\backslash \Omega_{i_3}$ and an index $i_4\in I$ such that $x_3\in W^s(\Omega_{i_4})\backslash \Omega_{i_4}$. This means that $\Omega_{i_1}\succ \Omega_{i_2}\succ \Omega_{i_3}\succ \Omega_{i_4}$. Because of the finiteness of the index set $I$, this procedure must be terminated in a finite step,  say $(m-1)$-step. Thus
$$
\Omega_{i_1}\succ\Omega_{i_2}\succ\cdots\succ \Omega_{i_m}
$$
holds, and either $i_m=i_1$ or $\Omega_{i_m}$ is Liapunov stable. An $(m-1)$-cycle is formed in the former case, and $\mu(\Omega_{i_j})=0$ for $j=1,2,\cdots, m-1$ is obtained by the topological transitivity of basic sets and Theorem \ref{saddle}.

It is well-known that $\Omega-$stablility is  equivalent to Axiom A and no-cycles condition (see \cite{Smale, Pugh, Hayashi}); Structural stability is equivalent to Axiom A and strong transversality condition, which implies no--cycles condition holds (see \cite{Smale, Robinson, Hayashi, Wen}); and separated star flow is equivalent to Axiom A and no--cycles condition (see \cite{WenGan, Wen}). Thus for each of these systems, limit measures live on Liapunov stable basic sets. This completes the proof.
\end{proof}
\begin{remark}
If we delete the assumption that there is no cycle in Corollary {\rm\ref{cor}}, then it holds  that any limit measure sits on Liapunov equilibria, or closed orbits or saddle cycles.
\end{remark}
In the following, we shall provide a series of examples whose concentrations of limit measures are precisely portrayed. Note that if we replace $b(x)$ by $\frac{b(x)}{\sqrt{1+|b(x)|^2}}$ then the two systems corresponding to (\ref{unpersys}) are topologically equivalent. Thus the modified drift terms $\frac{b(x)}{\sqrt{1+|b(x)|^2}}$ in all following constructed examples satisfy the existence conditions about drift terms  of FWULDP
 (see \cite[Theorem 3.1, p.135]{FW2}). If we assume the diffusion matrices also satisfies the corresponding conditions of \cite[Theorem 3.1, p.135]{FW2}, then {\rm(\ref{itodff})} admits FWULDP. We have also checked all following examples satisfy the existence conditions of FWULDP in the paper \cite[Theorem 2.1]{Yang}   as long as the order of the diffusion terms $\sigma$ tending to infinity are properly limited.

Let $H: \mathbb{R}^2\rightarrow \mathbb{R}$ be a $C^2$ Hamiltonian function, which determines a Hamiltonian system by the Hamiltonian vector field $\mathcal{H}(H):=(\frac{\partial H}{\partial x_2},-\frac{\partial H}{\partial x_1})^*$. Now we consider the system
\begin{equation}\label{PlanS}
\dot{x}=\mathcal{H}(H)-F(H)\nabla H
\end{equation}
where $F: {\rm Dom}(H)\rightarrow \mathbb{R}$ is a $C^1$ function.

\begin{prp}\label{EquiC}
Let $h_1\leq h_2$ and $F|_{[h_1,h_2]}\equiv 0$. Suppose that $H^{-1}\big([h_1,h_2]\big)$ is full of nontrivial periodic orbits of the Hamiltonian system determined by $H$. Then $H^{-1}\big([h_1,h_2]\big)$ is an equivalent class of the corresponding perturbed system {\rm(\ref{itodff})} of {\rm(\ref{PlanS})}.
\end{prp}
\begin{proof}
Let $x,y\in H^{-1}\big([h_1,h_2]\big)$ with $H(x)=a$ and $H(y)=b$. Then $a,b\in [h_1,h_2]$. If $a=b$, then $x$ and $y$ belong to the same periodic orbit $H^{-1}(a)$. Thus it is obvious that $x\thicksim y$.

Suppose that $a\neq b$, without loss of generality, we may assume that $a<b$. We shall prove that for any given $\eta$, there exist $T>0$ and $\varphi\in {\bf C}_{T}$ with $\varphi(0)=x$ and $\varphi(T)=y$ such that  $S_T(\varphi)<\eta$, which will imply that $V(x,y)=0$. Now we construct the following auxiliary system
\begin{equation}\label{AuxS}
\dot{x}=\mathcal{H}(H)-\lambda\nabla H
\end{equation}
where $\lambda<0$ with $|\lambda|$ sufficiently small. Let $\varphi^{\lambda}(t)$ be the solution of (\ref{AuxS}) passing through $x$. Then
\begin{equation}\label{Deriv}
\frac{d}{dt}H\big(\varphi^{\lambda}(t)\big)=-\lambda \big|\nabla H\big(\varphi^{\lambda}(t)\big)\big|^2.
\end{equation}
This shows that $H\big(\varphi^{\lambda}(t)\big)$ is increasing on $t\geq 0$. We claim that $\varphi^{\lambda}(t)$ will leave $H^{-1}\big([h_1,h_2]\big)$ from $H^{-1}\big(h_2\big)$ as $t$ is increasing. Otherwise, $\varphi^{\lambda}(t)\in H^{-1}\big([h_1,h_2]\big)$ for $t\geq 0$. The LaSalle invariant principle implies that the omega limit set of the orbit $\varphi^{\lambda}(t)$ must lie on $\{x\in H^{-1}\big([h_1,h_2]\big): \nabla H(x)=0\}$. However, by assumption, $H^{-1}\big([h_1,h_2]\big)$ is full of nontrivial periodic orbits, hence $\nabla H(x)\neq 0$ on $H^{-1}\big([h_1,h_2]\big)$, that is, $\{x\in H^{-1}\big([h_1,h_2]\big): \nabla H(x)=0\}$ is empty, a contraction. This proves the claim. Therefore, there exists a time $\tilde{t}>0$ such that $\varphi^{\lambda}(\tilde{t})\in H^{-1}\big(b\big)$. We extend the domain of definition of $\varphi^{\lambda}$ up
to $T$ as the solution of (\ref{PlanS}) with $\Psi_{T-\tilde{t}}\big(\varphi^{\lambda}(\tilde{t})\big)=y$,
which does not increase the value of $S$.

From the compactness of $H^{-1}\big([h_1,h_2]\big)$, $F|_{[h_1,h_2]}\equiv 0$ and (\ref{Deriv}), it follows that there is a constant $M>0$ such that
\begin{align*}
  S_T\big(\varphi^{\lambda}\big)= & S_{\tilde{t}}\big(\varphi^{\lambda}\big)\\
  \leq & M\int_0^{\tilde{t}}\big|\dot{\varphi}^{\lambda}-b(\varphi^{\lambda})\big|^2dt\\
  = & M\lambda^2\int_0^{\tilde{t}}\big|\nabla H\big(\varphi^{\lambda}(t)\big)\big|^2dt\\
   =& -M\lambda\int_0^{\tilde{t}}\frac{d}{dt}H\big(\varphi^{\lambda}(t)\big)dt\\
  = & -M\lambda(b-a)\\
  \leq & -M\lambda(h_2-h_1).
  \end{align*}
Thus, as $-\lambda<\frac{\eta}{M(h_2-h_1)}$, we have $S_T(\varphi^{\lambda})<\eta$. This proves that $V(x,y)=0$. Let $\lambda>0$ and repeat the above process, we can get that $V(y,x)=0$. This completes the proof.
\end{proof}

\begin{exmp}
Consider two dimensional SDEs:
\begin{equation}\label{ClosedO}
\left\{\begin{aligned}
dx_1=& \big(-x_2 + x_1f(x_1^2+x_2^2)\big)dt+\varepsilon \sigma_1(x_1,x_2)dw_t^1\\
dx_2=& \big(x_1+ x_2f(x_1^2+x_2^2)\big)dt+\varepsilon \sigma_2(x_1,x_2)dw_t^2,
\end{aligned}
\right.
\end{equation}
where
\[f(s)\begin{cases}
= 0,\quad & s\in [0,1]\cup\{4\},\\
> 0,\quad & s\in (1,4),\\
<0, \quad & s>4
\end{cases} \]
is a $C^2$ function.

By Proposition {\rm \ref{EquiC}} and the continuity of $V(x,y)$, the unit disk $D:=\{(x_1,x_2): r\leq 1\}$ is an equivalent class. It is easy to see that $B(\Psi)=D\bigcup\{r=\sqrt{x_1^2+x_2^2}=2\}$. Moreover, $D$ is a repeller for the corresponding flow $\Psi$, and the circle $r=2$ is an attractor (see Figure {\rm\ref{eg551}}).
\begin{figure}[h]
  \centering
  % Requires \usepackage{graphicx}
  \includegraphics[width=0.32\textwidth]{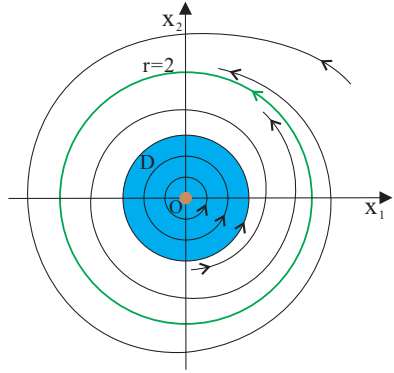}\\
  \caption{Phase Portrait for the noiseless system of (\ref{ClosedO}).}\label{eg551}
\end{figure}
Applying Proposition {\rm\ref{PR}} and Theorem {\rm\ref{repeller}}, we obtain that
$$\mu^{\varepsilon}\xlongrightarrow{w} \frac{1}{2\pi}\int_{0}^{2\pi}\delta_{\Psi_{t}\big((2,0)\big)}(\cdot)dt\ {\rm as }\ \varepsilon\rightarrow 0.$$

 Suppose that $C^2$ function $f$ in {\rm(\ref{ClosedO})} is replaced by
 \[f(s)\begin{cases}
< 0,\quad & s\in [0,1),\\
\equiv 0,\quad & s\in [1,4],\\
<0, \quad & s>4.
\end{cases} \]
Then it follows from Proposition {\rm\ref{EquiC}} that  the annulus region $D_1:=\{(x_1,x_2): 1\leq r\leq 2\}$ is is an equivalent class, which is asymptotically orbitally stable from the exterior and asymptotically orbitally unstable from the interior. Besides, the origin $O$ is an attractor and $B(\Psi)=D_1\cup\{O\}$ (see Figure {\rm\ref{eg552}}).
\begin{figure}[h]
  \centering
  % Requires \usepackage{graphicx}
  \includegraphics[width=0.32\textwidth]{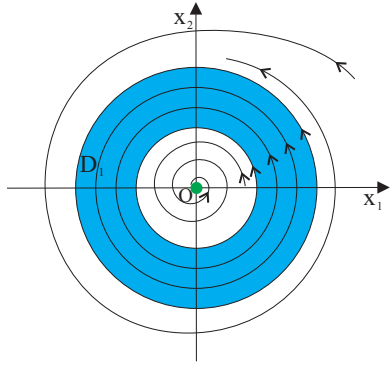}\\
  \caption{Phase Portrait for the noiseless system of (\ref{ClosedO}).}\label{eg552}
\end{figure}
Applying Proposition {\rm\ref{PR}} and Theorem {\rm\ref{saddle}}, we get that
 $$\mu^{\varepsilon}\xlongrightarrow{w} \delta_O(\cdot)\ {\rm as }\ \varepsilon\rightarrow 0.$$

\end{exmp}

\begin{exmp}\label{Ex3}
Consider two dimensional $C^2$ SDEs possessing infinite equivalent classes.
\begin{equation}\label{Accumulation}
\left\{\begin{aligned}
dx_1=& \Big(-x_2 + x_1\big(1-x_1^2-x_2^2\big)^5\sin^2 \big(1-x_1^2-x_2^2\big)^{-1}\Big)dt+\varepsilon \sigma_1(x_1,x_2)dw_t^1\\
dx_2=& \Big(x_1 + x_2\big(1-x_1^2-x_2^2\big)^5\sin^2 \big(1-x_1^2-x_2^2\big)^{-1}\Big)dt+\varepsilon \sigma_2(x_1,x_2)dw_t^2.
\end{aligned}
\right.
\end{equation}

\begin{prp}
For {\rm(\ref{Accumulation})}, we have
\begin{equation}\label{convergence}
\mu^{\varepsilon}\xlongrightarrow{w} \frac{1}{2\pi}\int_{0}^{2\pi}\delta_{\Psi_{t}\big((1,0)\big)}(\cdot)dt\ {\rm as }\ \varepsilon\rightarrow 0.
\end{equation}
The same result holds if the function $\big(1-x_1^2-x_2^2\big)^5\sin^2 \big(1-x_1^2-x_2^2\big)^{-1}$ in {\rm(\ref{Accumulation})} is replaced by
a $C^2$ function $f$ with the following properties
\[f(s)=\begin{cases}
 f^{-}(s)\geq 0,\quad & s\in [0,1),\\
 0,\quad & s=1,\\
f^{+}(s)\leq 0, \quad & s>1
\end{cases} \]
where
\[f^{-}(s)\begin{cases}
 \equiv 0,\quad & s\in I:=\bigcup_{n=1}^{\infty} [1-\frac{1}{(2n-1)\pi},1-\frac{1}{2n\pi}],\\
> 0, \quad & s\in [0,1)\backslash I;
\end{cases} \]
and
\[f^{+}(s)\begin{cases}
 \equiv 0,\quad & s\in J:=\bigcup_{n=1}^{\infty} [1+\frac{1}{2n\pi}, 1+\frac{1}{(2n-1)\pi}],\\
< 0, \quad &(1,+\infty)\backslash J.
\end{cases} \]
\end{prp}
\begin{proof}
Global dynamics of the corresponding determined system of {\rm(\ref{Accumulation})} is analyzed as follows. The origin $O(0,0)$ is a repeller. The unit circle $S$ is a closed orbit, which is accumulated by a series of semistable limit cycles:
$$\Gamma_n: r=\sqrt{1-\frac{1}{n\pi}},\  n =\pm 1, \pm 2, \pm 3, \cdots.$$
The Birkhoff center $B(\Psi)=\{O\}\bigcup_{n=-\infty}^{n=+\infty}\Gamma_n$ with $\Gamma_0=S$. The global phase portraits are sketched in Figure \ref{eg571}.
\begin{figure}[h]
  \centering
  % Requires \usepackage{graphicx}
  \includegraphics[width=0.35\textwidth]{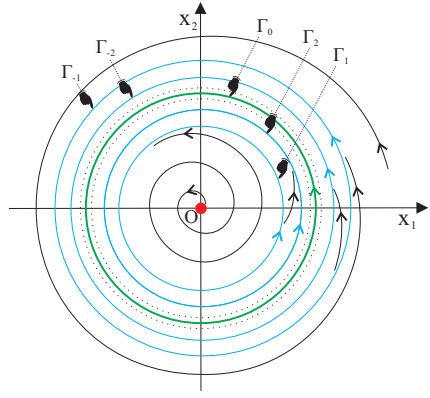}\\
  \caption{Phase Portrait for the noiseless system of (\ref{Accumulation}).}\label{eg571}
\end{figure}
Theorem \ref{repeller} applies to the repeller origin, we know that there is no concentration on the origin $O$. Let $R_n$ denote the annulus surrounded by $\Gamma_n$ and $\Gamma_{-n}$ for any positive integer $n$. Then $R_n$ is an attractor. Utilizing Theorem \ref{saddle} to semistable limit cycles $\{\Gamma_m: \  m =\pm 1, \pm 2, \pm 3, \cdots, \pm (n-1)\}$ and the attractor $R_n$, we conclude that there is no concentration on semistable limit cycles $\{\Gamma_m: \  m =\pm 1, \pm 2, \pm 3, \cdots, \pm (n-1)\}$. Since $n$ is arbitrary, we get (\ref{convergence}).

Moreover, let the  $C^2$ function $f$ satisfies the given properties and define the annuli
 $$\mathcal{A}_{-n}:=\left\{\sqrt{1-\big((2n-1)\pi\big)^{-1}}\leq r\leq \sqrt{1-\big(2n\pi\big)^{-1}} \right\},$$
 $$\mathcal{A}_{n}:=\left\{\sqrt{1+\big(2n\pi\big)^{-1}}\leq r\leq \sqrt{1+\big((2n-1)\pi\big)^{-1}} \right\},\ {\rm and}$$
 $$\mathcal{R}_{n}:=\left\{\sqrt{1-\big((2n-1)\pi\big)^{-1}}\leq r\leq \sqrt{1+\big((2n-1)\pi\big)^{-1}} \right\}$$
for $n=1,2,\cdots$. It follows that $B(\Psi)=\{O\}\bigcup_{n=-\infty}^{n=+\infty}\mathcal{A}_n$ with $\mathcal{A}_0=S$.
The global phase portraits are sketched in Figure \ref{eg572}.
\begin{figure}[h]
  \centering
  % Requires \usepackage{graphicx}
  \includegraphics[width=0.35\textwidth]{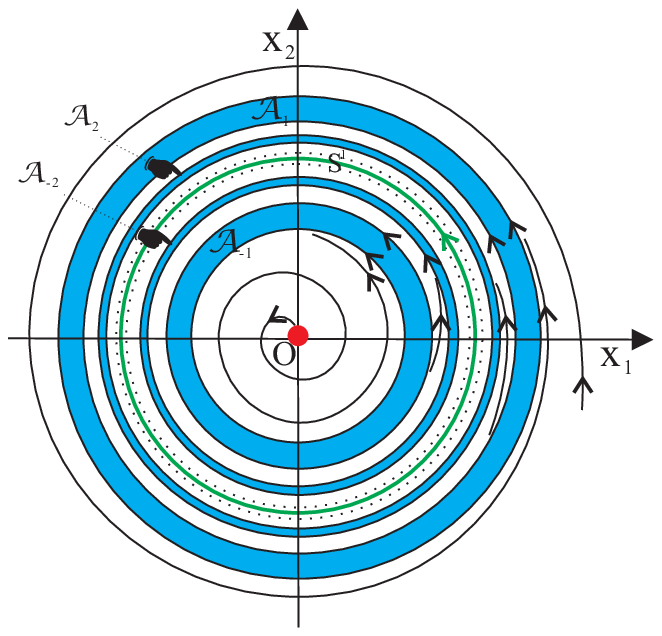}\\
  \caption{Phase Portrait for the noiseless system of (\ref{Accumulation}).}\label{eg572}
\end{figure}
Then by Proposition \ref{EquiC}, $\mathcal{A}_{-n}$ and $\mathcal{A}_{n}$ are equivalent classes for each $n$. Replacing by $\Gamma_n$, $\Gamma_{-n}$ and $R_n$ by $\mathcal{A}_{n}$, $\mathcal{A}_{-n}$ and $\mathcal{R}_{n}$, respectively, we obtain
 {\rm(\ref{convergence})} in the same manner.
\end{proof}
\end{exmp}

\begin{exmp}\label{Ex4}
Let $H(x_1,x_2):= \frac{x_2^2}{2}+\frac{x_1^4}{4}-\frac{x_1^2}{2}\in [-\frac{1}{4}, +\infty)$. Consider two dimensional nondegenerate SDEs:
\begin{equation}\label{FigureEight}
dx=[\mathcal{H}(H)-F(H)\nabla H]dt+\varepsilon \sigma(x)dw_t
\end{equation}
where $C^1$ function $F$ is taken
\[F_i(s)=\begin{cases}
(-1)^i|s|^3, \quad & s\in[-\frac{1}{4},0),\\
s^5\sin^2\frac{\pi}{s},\quad & s\in [0,1],\\
>0, \quad & s\in (1,2),\\
1,\quad & s\in [2,+\infty),\\
\end{cases} \]
for $i=1,2$, or
\[F_j(s)=\begin{cases}
(-1)^j|s|^3, \quad & s\in[-\frac{1}{4},0),\\
G(s),\\
\end{cases} \]
for $j=3,4$. Here $G:[0,+\infty)\rightarrow [0,1]$ is a $C^{\infty}$ function satisfying that $G(0)=0$ and there exist a sequence of intervals $I_n=[a_n,b_n]\subset (0,1), n=1,2,3,\cdots$ with $I_{n+1}\subset I_n$, $\lim_{n\rightarrow \infty}b_n=0$ and $b_1=1$ such that $G(s)=0$ for any $s\in I=\bigcup_{n=1}^{\infty}I_n$ and $G(s)>0$ for any $s\in (0,+\infty)\backslash I$. Note that such a $G$ can be constructed by cut-off functions.
\end{exmp}

For these systems, we have the following conclusions.
\begin{prp}
For $i=1,3$, we have
\begin{equation}\label{convergence1}
\mu^{\varepsilon}\xlongrightarrow{w} \delta_O(\cdot)\ {\rm as }\ \varepsilon\rightarrow 0;
\end{equation}
and for $i=2,4$, we have
\begin{equation}\label{convergence2}
\mu^{\varepsilon}\xlongrightarrow{w} \lambda_1\delta_{(1,0)}(\cdot)+\lambda_2\delta_{(-1,0)}(\cdot)\ {\rm as }\ \varepsilon\rightarrow 0
\end{equation}
with $\lambda_1+\lambda_2=1$.
\end{prp}
\begin{proof}
For $i=1$, by definition of solution, $H^{-1}\big(n^{-1}\big)$ is a limit cycle of {\rm(\ref{PlanS})} for $n=1,2,3,\cdots$, which is asymptotically orbitally stable from the exterior and asymptotically orbitally unstable from the interior.  These limit cycles accumulate on the homoclinic cycle, a figure-eight curve. Using $H$ as a Liapunov function, we can prove that any limit set of $\Psi$ lies on the zero set of $F(H)|\nabla H|^2$ (see Figure \ref{eg58i1}).
\begin{figure}[h]
  \centering
  % Requires \usepackage{graphicx}
  \includegraphics[width=0.35\textwidth]{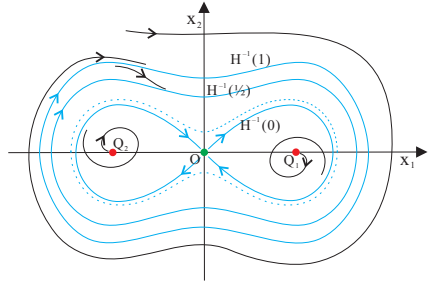}\\
  \caption{Phase Portrait for the noiseless system of (\ref{FigureEight}) as $i=1$.}\label{eg58i1}
\end{figure} Thus $B(\Psi)=\bigcup_{n=1}^{+\infty}H^{-1}\big(n^{-1}\big)\bigcup H^{-1}(0)\bigcup\{(\pm 1,0)\}$.

 Let $\mu^{\varepsilon}_k\xlongrightarrow{w} \mu$ as $\varepsilon_k\rightarrow 0$. Using $H+\frac{1}{4}$ as a Liapunov function, we can verify that $(\pm 1,0)$ is a repeller. Theorem \ref{repeller} implies that there is no concentration of $\mu$ on the equilibria $(\pm 1,0)$. Again using the Liapunov method, we can prove that $H^{-1}\big([-4^{-1},n^{-1}]\big)$ is an attractor for each $n$. Utilizing Theorem \ref{saddle} to semistable limit cycles $\{H^{-1}\big(m^{-1}\big): \  m = 1, 2, 3, \cdots, n-1\}$ and the attractor $H^{-1}\big([-4^{-1},n^{-1}]\big)$, we conclude that there is no concentration of $\mu$ on semistable limit cycles $\{H^{-1}\big(m^{-1}\big): \  m = 1, 2, 3, \cdots, n-1\}$. Since $n$ is arbitrary, we obtain that ${\rm supp}(\mu)\subset H^{-1}(0)$. Since $H^{-1}(0)$ consists of two homoclinic orbits connecting the origin $O$, it follows from the invariance of $\mu$ with respect to $\Psi$ that $\mu$ must be $\delta_O(\cdot)$. This proves (\ref{convergence1}).

For $i=3$, the annular region $H^{-1}\big([a_n,b_n]\big)$ is full of nontrivial periodic orbits of {\rm(\ref{PlanS})} for $n=1,2,3,\cdots$, which is asymptotically orbitally stable from the exterior and asymptotically orbitally unstable from the interior.  These annular regions accumulate on the homoclinic cycle. Similarly, $B(\Psi)=\bigcup_{n=1}^{+\infty}H^{-1}\big([a_n,b_n]\big)\bigcup H^{-1}(0)\bigcup\{(\pm 1,0)\}$.

 Let $\mu^{\varepsilon}_k\xlongrightarrow{w} \mu$ as $\varepsilon_k\rightarrow 0$. In the same manner we can prove that there is no concentration of $\mu$ on the equilibria $(\pm 1,0)$ and that $H^{-1}\big([-4^{-1},b_n]\big)$ is an attractor for each $n$. By Proposition \ref{EquiC},   the annular region $H^{-1}\big([a_n,b_n]\big)$ is an equivalent class for each $n$. Applying Theorem \ref{saddle} to the annular regions
$\{H^{-1}\big([a_m,b_m]\big): \  m = 1, 2, 3, \cdots, n-1\}$ and the attractor $H^{-1}\big([-4^{-1},b_n]\big)$, we conclude that there is no concentration of $\mu$ on the annular regions $\{H^{-1}\big([a_m,b_m]\big): \  m = 1, 2, 3, \cdots, n-1\}$. Since $n$ is arbitrary, we obtain that ${\rm supp}(\mu)\subset H^{-1}(0)$ and (\ref{convergence1}) holds.

Suppose that $i=2,4$. Then the same procedure proves that there is no concentration of $\mu$ on  both $H^{-1}\big(n^{-1}\big)$ and $H^{-1}\big([a_n,b_n]\big)$ for each $n$. Again using $H+\frac{1}{4}$ as a Liapunov function, we can verify that $(\pm 1,0)$ are attractors for both cases. Since every trajectory  of $\Psi$ in the interior of $H^{-1}(0)$ other than $(\pm 1,0)$ converges to $H^{-1}(0)$ as $t\rightarrow -\infty$ and to $(\pm 1,0)$ as $t\rightarrow +\infty$ (see Figure \ref{eg58i2}).
\begin{figure}[h]
  \centering
  % Requires \usepackage{graphicx}
  \includegraphics[width=0.35\textwidth]{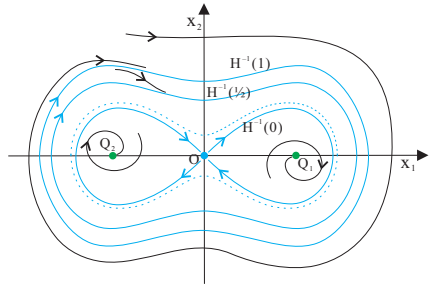}\\
  \caption{Phase Portrait for the noiseless system of (\ref{FigureEight}) as $i=2$.}\label{eg58i2}
\end{figure} Applying Theorem \ref{saddle} to $H^{-1}(0)$ and $(\pm 1,0)$, we get that $\mu\big(H^{-1}(0)\big)=0$. Finally, ${\rm supp}(\mu)=\{(\pm 1,0)\}$, which proves (\ref{convergence2}).
\end{proof}

Examples \ref{Ex3} and \ref{Ex4} show that there is no concentration on semistable limit cycles. Actually, this result holds for any planar systems by Theorem\ref{saddle}. However, the following two examples illustrates limit measure can be supported on saddle-node (semistable equilibrium). The essential reason is that there is a cycle connecting it.

\begin{exmp}
Consider two dimensional SDEs:
\begin{equation}\label{Saddle-node1}
\left\{\begin{aligned}
dx_1=& \big[x_1(1-x_1^2-x_2^2)-x_2(1+x_1)\big]dt+\varepsilon \sigma_1(x_1,x_2)dw_t^1\\
dx_2=& \big[x_1(1+x_1)+x_2(1-x_1^2-x_2^2)\big]dt+\varepsilon \sigma_2(x_1,x_2)dw_t^2.
\end{aligned}
\right.
\end{equation}
The unperturbed system has equilibria $O(0,0)$ and $Q(-1,0)$, which make up the Birkhoff center $B(\Psi)$.  $O$ is a repeller and $Q$ is a saddle-node which is  connected by a homoclinic orbit on the unit circle, see Figure {\rm\ref{eg510}}.
\begin{figure}[h]
  \centering
  % Requires \usepackage{graphicx}
  \includegraphics[width=0.35\textwidth]{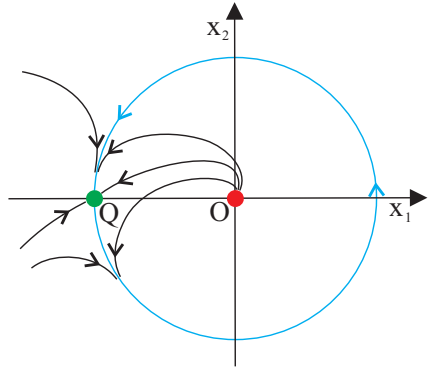}\\
  \caption{Phase Portrait for the noiseless system of (\ref{Saddle-node1}).}\label{eg510}
\end{figure}
By Theorem {\rm\ref{repeller}}, there is no concentration on the origin $O$. Thus
$$\mu^{\varepsilon}\xlongrightarrow{w} \delta_Q(\cdot)\ {\rm as }\ \varepsilon\rightarrow 0.$$

\end{exmp}

\begin{exmp}\label{Saddle-nodeC}
Consider two dimensional SDEs:
\begin{equation}\label{Saddle-node}
\left\{\begin{aligned}
dx_1=& \big[x_1\big(x_1^2+x_2^2\big)\big(1-x_1^2-x_2^2\big)-2x_2^3\big]dt+\varepsilon \sigma_1(x_1,x_2)dw_t^1\\
dx_2=& \big[x_2\big(x_1^2+x_2^2\big)\big(1-x_1^2-x_2^2\big)+2x_1x_2^2\big]dt+\varepsilon \sigma_2(x_1,x_2)dw_t^2.
\end{aligned}
\right.
\end{equation}
The corresponding determined system has equilibria $O(0,0)$ and $(\pm 1,0)$, which comprise the Birkhoff center $B(\Psi)$.  $O$ is a repeller and $(\pm 1,0)$ are saddle-nodes which are connected by two entire orbits on the unit circle,
see Figure {\rm\ref{eg511}}.
\begin{figure}[h]
  \centering
  % Requires \usepackage{graphicx}
  \includegraphics[width=0.33\textwidth]{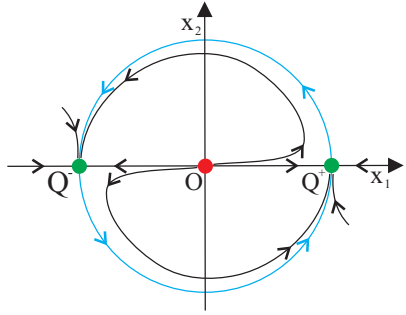}\\
  \caption{Phase Portrait for the noiseless system of (\ref{Saddle-node}).}\label{eg511}
\end{figure}
By Theorem {\rm\ref{repeller}}, there is no concentration on the origin $O$. Thus {\rm(\ref{convergence2})} holds.
\end{exmp}

\begin{exmp}\label{Kifer}
Consider two dimensional determined system:
\begin{equation}\label{KiferE}
\left\{\begin{aligned}
\frac{dx_1}{dt}=& x_1\left(x_1^2+x_2^2\right)\left(1-x_1^2-x_2^2\right)-x_1x_2\left(\sqrt{x_1^2+x_2^2}-x_1\right)\\
\frac{dx_2}{dt}=& x_2\left(x_1^2+x_2^2\right)\left(1-x_1^2-x_2^2\right)+x_1^2\left(\sqrt{x_1^2+x_2^2}-x_1\right).
\end{aligned}
\right.
\end{equation}
The system {\rm(\ref{KiferE})} has equilibria  $A(1,0)$, $B(0,1)$, $C(0,-1)$ and $O(0,0)$, $B$ is a stable node, $C$ is a saddle,  $A$ is a saddle-node and $O$ is a repeller, $S^1$ is invariant and attracts all points except the origin $O$,
see Figure {\rm\ref{egKifer}}.
\begin{figure}[h]
  \centering
  % Requires \usepackage{graphicx}
  \includegraphics[width=0.32\textwidth]{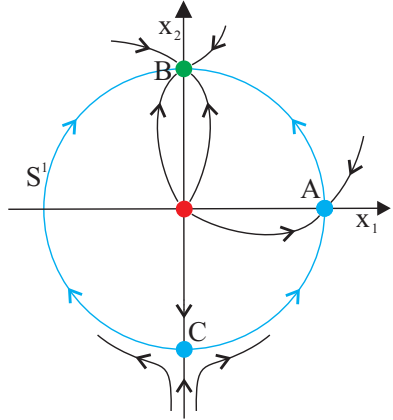}\\
  \caption{Phase Portrait for the system (\ref{KiferE}).}\label{egKifer}
\end{figure}
When the flow is restricted on $S^1$, it is topologically equivalent to that given by Kifer {\rm\cite [Remark 5.4, p.90-91]{Kifer}}. However, he chose a perturbation such that its limit measure is supported at the saddle-node $A$. But when {\rm(\ref{KiferE})} is perturbed by any nondegenerate while noise, its limit measure will sit on the stable node $B$ by Theorems {\rm\ref{repeller}} and {\rm\ref{saddle}}, which is different from Kifer's.
\end{exmp}

\begin{exmp}\label{van der Pole}
Consider the stochastic perturbation of the variational equation of the forced van der Pol oscillator:
\begin{equation}\label{Saddle-nodeB}
\left\{\begin{aligned}
du=& \big[u-\sigma v-u\big(u^2+v^2\big)\big]dt+\varepsilon \sigma_1(u,v)dw_t^1\\
dv=& \big[\sigma u+v-v\big(u^2+v^2\big)-\gamma\big]dt+\varepsilon \sigma_2(u,v)dw_t^2.
\end{aligned}
\right.
\end{equation}
The corresponding determined system of {\rm(\ref{Saddle-nodeB})} is the well-known variational equation of the forced van der Pol oscillator. In a neighborhood of the parameters $(\sigma,\gamma)=(\frac{1}{2},\frac{1}{2})$, Guchenheimer and Holmes {\rm\cite [p.70-74]{GH}} completely summarized the saddle-node bifurcations and Hopf bifurcations and drew $23$ phase portrait diagrams of the corresponding determined system in Figure {\rm2.1.3} of {\rm\cite [p.72]{GH}}.
\end{exmp}

Combing Theorems \ref{repeller} and \ref{saddle} with Figure 2.1.3 of {\rm\cite [p.72]{GH}}, we conclude the following.
\begin{prp}\label{Van der Pole}

{\rm (i)} There are $8$ diagrams having sources on which there is no concentration by Theorem {\rm\ref{repeller}};

{\rm (ii)} there are $7$ diagrams having saddles on which there is no concentration by Theorem {\rm\ref{saddle}} {\rm(}because there exist entire orbits connecting them to sinks or stable limit cycles{\rm)};

{\rm (iii)} there are $5$ diagrams having saddle-nodes on which there is no concentration by Theorem {\rm\ref{saddle}} {\rm(}because there exist entire orbits connecting them to sinks or stable limit cycles{\rm)};

{\rm (iv)} there is one diagram having a codimension $2$ equilibrium  on which there is no concentration by Theorem {\rm\ref{saddle}} {\rm(}because there exists a entire orbit connecting it to a sink{\rm)};

{\rm (v)} there are $2$ diagrams having degenerate equilibria with a homoclinic orbit, which were sketched in the second line from bottom of {\rm\cite [Figure 2.1.3, p.72]{GH}}, limit measure may be supported on them as in the Example {\rm\ref{Saddle-nodeC}}.
\end{prp}

\begin{exmp}
Consider three dimensional SDEs:
\begin{equation}\label{Cooperative}
dX^{\varepsilon}_t=b(X^{\varepsilon}_t)dt+\varepsilon \sigma(X^{\varepsilon}_t)dw_t
\end{equation}
where $\sigma(X)$ is nondegenerate, $\sigma(-X)=\sigma(X)$ and $b(-X)=-b(X)$ for any $X\in \mathbb{R}^3$, $b(X)=b(x_1,x_2,x_3)$ is a cooperative and irreducible vector field given by
\begin{equation}\label{CI}
\dot{X}=b(x_1,x_2,x_3)=\left(-x_1+\frac{x_3}{1+|x_3|}, x_1-x_2, x_2-\frac{1}{2}x_3\right)^*.
\end{equation}
The equilibria set of {\rm(\ref{CI})} is $O(0,0,0)$, $P_{\pm}=\pm(1,1,1)$, $O$ is a saddle and $P_{\pm}$ are asymptotically stable. There is a connecting orbit from $O$ to $P_+(P_-)$. It follows from {\rm\cite {Selgrade}} that every trajectory of {\rm(\ref{CI})} tends to an equilibrium, see Figure \ref{seg515}.
\begin{figure}[h]
  \centering
  % Requires \usepackage{graphicx}
  \includegraphics[width=0.35\textwidth]{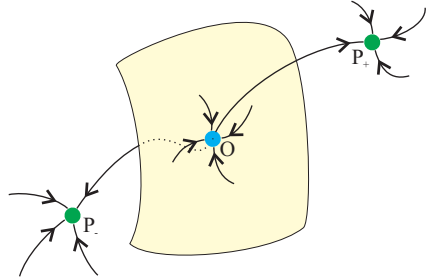}\\
  \caption{Phase Portrait for the vector field (\ref{CI}).}\label{seg515}
\end{figure}
By symmetry of $b, \sigma$,  the normalized solution $p^{\varepsilon}$  of the Fokker--Planck equation corresponding to {\rm(\ref{Cooperative})} satisfies $p^{\varepsilon}(-X)=p^{\varepsilon}(X)$ for any $X\in \mathbb{R}^3$. Suppose that $\mu^{\varepsilon}$ is the stationary measure of {\rm(\ref{Cooperative})}. Then  $p^{\varepsilon}$ is its density function.  Let $\mu^{\varepsilon}\xlongrightarrow{w}\mu$ as $\varepsilon\rightarrow 0$. Then by Theorem {\rm\ref{saddle}}, $\mu(O)=0$. It follows from the symmetry of $p^{\varepsilon}$ that as $\varepsilon \rightarrow 0$,
$$\mu^{\varepsilon}\xlongrightarrow{w}\frac{1}{2}\delta_{P_-}+\frac{1}{2}\delta_{P_+}.$$
\end{exmp}

\begin{exmp}
Let $H=\frac{x_1^4}{4}-\frac{x_1^2}{2}+\frac{x_2^2}{2}$, $\mathcal{H}(H)=(x_2, -(x_1^3-x_1))^*$, $\nabla H=((x_1^3-x_1),x_2)^*$ and $b(X)=\mathcal{H}(H)-(H+\frac{1}{8})\nabla H$.
Consider two dimensional SDEs:
\begin{equation}\label{TwoCO}
dX^{\varepsilon}_t=b(X^{\varepsilon}_t)dt+\varepsilon \sigma(X^{\varepsilon}_t)dw_t
\end{equation}
where $\sigma(X)$ is nondegenerate, $\sigma(-X)=\sigma(X)$ for any $X\in \mathbb{R}^2$.  It is easy to see $b(-X)=-b(X)$ for any $X\in \mathbb{R}^2$. The corresponding perturbed system of {\rm(\ref{TwoCO})} is
\begin{equation}\label{TwoCO1}
\dot{X}=b(X).
\end{equation}
The equilibria set of {\rm(\ref{TwoCO1})} is $O(0,0)$, $P_{\pm}=(\pm1,0)$, $O$ is a saddle without any homoclinic orbit and $P_{\pm}$ are repellers. Besides, $H+\frac{1}{8}=0$ consists of two limit cycles $\Gamma_+$ and $\Gamma_-$, which are attractors, $B(\Psi)=\Gamma_{\pm}\cup \{O,P_{\pm}\}$, see Figure {\rm\ref{seg516}}.
\begin{figure}[h]
  \centering
  % Requires \usepackage{graphicx}
  \includegraphics[width=0.38\textwidth]{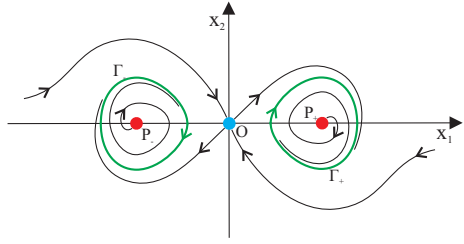}\\
  \caption{Phase Portrait for the noiseless system of (\ref{TwoCO}).}\label{seg516}
\end{figure}
By symmetry of $b, \sigma$,  the normalized solution $p^{\varepsilon}$  of the Fokker--Planck equation corresponding to {\rm(\ref{TwoCO})} satisfies $p^{\varepsilon}(-X)=p^{\varepsilon}(X)$ for any $X\in \mathbb{R}^2$. Suppose that $\mu^{\varepsilon}$ is the stationary measure of {\rm(\ref{TwoCO})}. Then  $p^{\varepsilon}$ is its density function.  Let $\mu^{\varepsilon}\xlongrightarrow{w}\mu$ as $\varepsilon\rightarrow 0$. Then by Theorems \ref{repeller} and \ref{saddle}, $\mu(P_{\pm})=\mu(O)=0$. It follows from the symmetry of $p^{\varepsilon}$ that
$$\mu^{\varepsilon}\xlongrightarrow{w} \frac{1}{2}\left\{\frac{1}{T}\int_{0}^{T}\delta_{\Psi_{t}\big((\sqrt{2},0)\big)}(\cdot)dt+\frac{1}{T}\int_{0}^{T}\delta_{\Psi_{t}\big((-\sqrt{2},0)\big)}(\cdot)dt\right\}\ {\rm as }\ \varepsilon\rightarrow 0.$$
\end{exmp}

The last example shows that Theorem \ref{repeller} is not valid if ${\bf\big(P_{\mathcal{R}}\big)}$ is violated.
\begin{exmp}\label{counterexample}
Consider stochastic differential equations
\begin{equation}\label{PRnot}
  \left\{
    \begin{array}{ll}
     \mathrm dx=[y-x(x^2+y^2-1)(x^2+y^2-2)]\mathrm dt+\varepsilon\mathrm dw^1_t,\\
      \mathrm dy=[-x-y(x^2+y^2-1)(x^2+y^2-2)]\mathrm dt+\varepsilon\mathrm dw^2_t.
    \end{array}
  \right.
\end{equation}
The phase portrait for the unperturbed system are shown in Figure {\rm\ref{seg517}} below.
\begin{figure}[h]
  \centering
  % Requires \usepackage{graphicx}
  \includegraphics[width=0.33\textwidth]{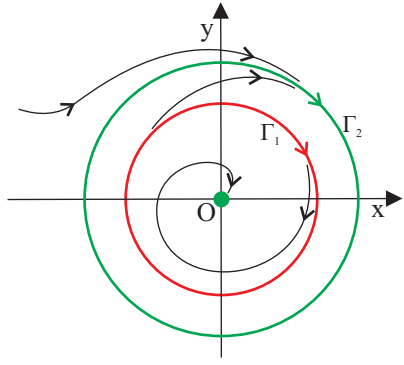}\\
  \caption{Phase Portrait for the noiseless system of (\ref{PRnot}).}\label{seg517}
\end{figure}
Here $O=(0,0)$ is a stable focus, $\Gamma_1=\{(x,y)\in\mathbb{R}^2:x^2+y^2=1\}$ is an unstable limit cycle
and $\Gamma_2=\{(x,y)\in\mathbb{R}^2:x^2+y^2=2\}$ is a stable limit cycle.
Set $\mathcal{R}=\{(x,y)\in\mathbb{R}^2:x^2+y^2\leq1\}$. Then $\mathcal{R}$ is a repeller for the unperturbed system.
Using the method in {\rm\cite{CJX}}, we can verify
that the unique invariant probability density for {\rm{\rm(\ref{PRnot})} } is
\[    p^\varepsilon(x,y)=C(\varepsilon)\exp\big(-2\varepsilon^{-2}U(x,y)\big),        \]
where $U(x,y)=(x^2+y^2)^3/6-3(x^2+y^2)^2/4+x^2+y^2$ and $C(\varepsilon)=\left(\int_{\mathbb{R}^2}\exp\big(-2\varepsilon^{-2}U(x,y)\big)\mathrm dx\mathrm dy\right)^{-1}<+\infty$.
Since $U(O)=0<\frac{1}{3}=U(x,y)$ for $(x,y)\in \Gamma_2$, by Laplace's method,
$\mu^{\varepsilon}\xlongrightarrow{w} \delta_{O}$ as $\varepsilon\rightarrow 0$. This implies that the limit measure is concentrated on the repeller $\mathcal{R}$.
\end{exmp}

\section{Appendix}

In this section, we first prove Proposition \ref{Hypothesis} under ${\bf(0_w)}$ and dissipativity condition of {\rm (\ref{unpersys}).
 The central objective of this appendix is to prove Theorem \ref{repeller} and Theorem \ref{mainatt} under FWULDP. The next proposition can be deduced easily from the compactness assumption in ${\bf(0_w)}$.
\begin{prp}\label{attain_1}
Let ${\bf(0_w)}$ hold. If $F\subset{\bf C}_{T} $  is  closed with $\cup_{t\in [0,T]}F^t:=\cup_{t\in [0,T]}\{\varphi(t):\varphi\in F \}$ being bounded, then there is a $\psi\in F$ such that
$$S_{T}(\psi)=\inf_{\varphi\in F}S_T(\varphi)=:S_T(F).$$
\end{prp}

\begin{prp}\label{invr_set}
 Suppose that {\rm (\ref{unpersys})} is dissipative. Then for any bounded set $A$, there is a positively invariant compact set $K\supset A$.
\end{prp}
\begin{proof}
Since $\Psi$ is dissipative, it has a global attractor $\mathcal{A}$ whose basin is $\mathbb{R}^r$.  Then it follows from Wilson
\cite [Theorem 3.2]{Wilson} that $\mathcal{A}$  has a $C^{\infty}$ Liapunov function $V:\mathbb{R}^r\rightarrow \mathbb{R}_+$ satisfying
$$V(\mathcal{A})=0,\ V (x)\wedge \langle -b(x),\nabla V(x)\rangle >0,\ x\in \mathbb{R}^r\setminus \mathcal{A},\ {\rm and} \lim_{|x|\rightarrow +\infty}V(x)=+\infty.$$
For any $\delta>0$, $V^{-1}([0,\delta])$ is a positively invariant compact set. By the compactness of $\overline{A}$ and $\lim_{|x|\rightarrow +\infty}V(x)=+\infty$, there exists an integer $m$ such that $\overline{A}\subset V^{-1}([0,m])=:K$, which is positively invariant and compact.
\end{proof}

{\bf Proof of Proposition \ref{Hypothesis}}.
We shall prove the proposition by contradiction. Suppose that $\{\varphi_n\}$ is a sequence of functions in $F$ satisfying $\lim_{n\rightarrow\infty}S_T(\varphi_n)=0$. Then by Proposition \ref{invr_set} we may choose a positively invariant compact set $K_1:=V^{-1}([0,m])\supset F^0$. According to the proof of Proposition \ref{invr_set}, $K_2:=V^{-1}([0,(m+1)])$ satisfies ${\rm Int} K_2\supset K_1$. Now for every $\varphi_n$, define $\phi_n$ to be $\varphi_n$ if $\varphi_n([0,T])\subset {\rm Int}K_2$. Otherwise, set $T_*^n=T_*(\varphi_n):=\inf\{t\in[0,T]:\varphi_n(t)\in \partial K_2\}$. Let $\phi_n$ be the same as $\varphi_n$ for $t\in [0,T_*^n]$ and extend it by the solution of system {\rm (\ref{unpersys})} for $t\in [T_*^n,T]$. We emphasize that in the second case, $\phi_n$ can not be a solution of system {\rm (\ref{unpersys})} since $K_1$ is positively invariant with respect to $\Psi$ and $\phi_n(0)=\varphi_n(0)\in K_1$. We also have $\phi_n([0,T])\subset K_2$ and $\lim_{n\rightarrow\infty}S_T(\phi_n)=0$. This means that after dropping finite terms, we may assume $\{\phi_n\}\subset \{\varphi\in{\bf C}_{T}:\varphi(t)\in K_2,t\in [0,T],S_{T}(\varphi)\leq 1\}=:\mathcal{F}$. Recall that $\mathcal{F}$ is compact in $\mathbf{C}_T$ by ${(\bf 0_w)}$. So we may assume again there is a $\phi_*\in \mathcal{F}$ such that $\lim_{n\rightarrow \infty}\phi_n=\phi_*$ with $\phi_*(0)\in K_1$. Indeed, $\phi_*$ can not be a solution of system {\rm (\ref{unpersys})}. Otherwise for sufficiently large $n$, we have $\phi_n$ close enough to $\phi_*$ so that $\phi_n([0,T])\subset {\rm Int}K_2$. However, from the construction of $\phi_n$ we know that this means $\phi_n=\varphi_n$. As a result, $\phi_*\in F$ by the closedness of $F$, which contradicts the hypothesis of this proposition. Finally, from the lower semi-continuity of $S_T$  and Remark \ref{Rem1}, we conclude that $\liminf_{n \rightarrow \infty} S_T(\phi_n)\geq S_T(\phi_*)>0$, a contradiction.

We now proceed to prove Theorem \ref{repeller} under FWULDP.

{\bf Proof of Theorem \ref{repeller} Under {\bf FWULDP}}. Choose $\delta>0$ such that $(\mathcal{R})_{\delta}$ is a fundamental neighborhood of $\mathcal{R}$.
By Lemma \ref{klofR}, there exists $\delta_1,\ \delta^*_2\in(0,\delta),\ \delta_1>\delta^*_2$
and $s_0>0$ such that $V(\partial (\mathcal{R})_{\delta_1},\partial (\mathcal{R})_{\delta^*_2})\geq s_0$. Set $\delta_2=\delta^*_2/2$. Applying Proposition \ref{reppp} to $(\mathcal{R})_{\delta}$ and $\overline{(\mathcal{R})_{\delta_1}}\backslash(\mathcal{R})_{\delta_2}$, we have
$T_0=\inf\big\{u\geq 0:\Psi_{t}\big(\overline{(\mathcal{R})_{\delta_1}}\backslash(\mathcal{R})_{\delta_2}\big)
\subset((\mathcal{R})_{\delta})^c,t\geq u\big\}<+\infty$. Let $F_0=\{\varphi\in {\bf C}_{T_0}:\varphi(0)\in \overline{(\mathcal{R})_{\delta_1}}\backslash(\mathcal{R})_{\delta_2}, \varphi(T_0)\in \overline{(\mathcal{R})_{\delta_1}}\}$ be a closed subset of ${\bf C}_{T_0}$. Then $F_0^0$ is bounded and $F_0$ does not contain any solution of system (\ref{unpersys}) by the definition of $T_0$.
Then by Hypothesis \ref{closeinfH}, we have $s_1=S_{T_0}(F_0)>0$.

We now follow the line of proof for Theorem \ref{repeller} under DZULDP to prove the same result under FWULDP. The central tasks are to prove that (\ref{repconout1}), (\ref{repconout2}) and (\ref{repconout3}). The same procedure proves (\ref{Esti1}).

For any $y\in \partial(\mathcal{R})_{\delta_1}$, let $F_{1,y}=\{\varphi\in {\bf C}_T^y: \varphi(t) \in \partial(\mathcal{R})_{\delta_2}
\:\textrm{for some} \ t\in[0,T]\}$ where $T$ is the same in the proof of Theorem \ref{repeller} under  DZULDP. For any $\tilde{\varphi}\in \{\varphi \in \mathbf{C}^y_T:\rho_T(\varphi,F_{1,y})<\delta^*_2-\delta_2\}$, there exists $\hat{\varphi}\in F_{1,y},\hat{t}\in[0,T]$ such that $\tilde{\varphi}(\hat{t})\in B(\hat{\varphi}(\hat{t}),\delta^*_2-\delta_2)\subset (\mathcal{R})_{\delta^*_2}$. From the continuity of $\tilde{\varphi}$, we know that there is some $\tilde{t}=\tilde{t}(\tilde{\varphi})\in [0,T]$ such that $\tilde{\varphi}(\tilde{t})\in \partial(\mathcal{R})_{\delta^*_2}$. By the definition of $V(\partial(\mathcal{R})_{\delta_1},\partial(\mathcal{R})_{\delta^*_2})$
and Lemma \ref{klofR}, we have $$\inf_{y\in \partial(\mathcal{R})_{\delta_1}}S_{T}^{y}\big(\{\varphi \in \mathbf{C}^y_T:\rho_T(\varphi,F_{1,y})<\delta^*_2-\delta_2\}\big)
\geq V(\partial(\mathcal{R})_{\delta_1},\partial(\mathcal{R})_{\delta^*_2})\geq s_0.$$
This proves that
$$\mathbb{F}_T^y(0.9s_0)\cap\{\varphi \in \mathbf{C}^y_T:\rho_T(\varphi,F_{1,y})<\delta^*_2-\delta_2\}=\emptyset$$
and
$$F_{1,y}\subset \{\varphi \in \mathbf{C}^y_T:\rho_T(\varphi,\mathbb{F}_T^y(0.9s_0))\ge\delta^*_2-\delta_2\}$$
 which means that $\{X^{\varepsilon}\in F_{1,y}\}\subset \{X^{\varepsilon}(0)=y,\rho_{T}(X^{\varepsilon},\mathbb{F}_{T}^{y}(0.9s_0))\geq\delta^*_2-\delta_2\},\forall y\in \partial(\mathcal{R})_{\delta_1}.$
By (\ref{ufldpupbb}) of ${\bf(II_u)}$,
there exists $\varepsilon_1>0$ such that for any $\varepsilon\in(0,\varepsilon_1)$
and $y\in \partial(\mathcal{R})_{\delta_1}$, we get
\begin{align}\label{Esti2_dis}
 \mathbb{P}_{y}(\tau\leq T)=\mathbb{P}_{y}(X^{\varepsilon}\in F_{1,y})
  &\leq\mathbb{P}_{y}\{\rho_{T}(X^{\varepsilon},\mathbb{F}_{T}^{y}(0.9 s_0))\geq\delta^*_2-\delta_2\} \nonumber\\
  &\leq\exp\{-\varepsilon^{-2}(0.9 s_0-0.1 s_0)\}\leq \exp\{-0.8s_0/\varepsilon^2\}.
\end{align}
Therefore, by (\ref{Esti1}) and (\ref{Esti2_dis}), for any $\varepsilon\in(0,\varepsilon_1)$,
\[\sup_{z\in \big(\overline{(\mathcal{R})_{\delta_1}}\big)^c}\mathbb{P}_z(X^{\varepsilon}_T \in (\mathcal{R})_{\delta_2} )\leq\sup_{y\in\partial(\mathcal{R})_{\delta_1}}\mathbb{P}_{y}(\tau\leq T)\leq \exp\{-0.8s_0/\varepsilon^2\}.
\]
Thus, the inequality (\ref{repconout1}) follows.

For any $z\in \overline{(\mathcal{R})_{\delta_1}}\backslash (\mathcal{R})_{\delta_2}$, let $F_{2,z}=\big\{\varphi\in {\bf C}_T^z: \varphi(T) \in \overline{(\mathcal{R})_{\delta_2}}\big\}$.
We claim that
\begin{equation}\label{Esti3_dis}
S_{T}\big(\{\varphi\in\mathbf{C}^z_T:\rho_T(\varphi,F_{2,z})<\delta^*_2-\delta_2 \}\big)\geq s_0\wedge s_1.
\end{equation}
Indeed, let $\tilde{\varphi}$ be any element in $\{\varphi\in\mathbf{C}^z_T:\rho_T(\varphi,F_{2,z})< \delta^*_2-\delta_2 \}$. Recall the definition of $F_{2,z}$, we conclude that $\tilde{\varphi}(T)\in \overline{(\mathcal{R})_{\delta^*_2}}$.
If $S_{T}(\tilde{\varphi})<s_1$, then by the definition of $s_1$,
we have $\tilde{\varphi}(T_0)\notin \overline{(\mathcal{R})_{\delta_1}}$.
Thanks to the continuity of $\tilde{\varphi}$ and
$\tilde{\varphi}(T)\in \overline{(\mathcal{R})_{\delta^*_2}}$,
there exist $t_1,t_2\in(T_0,T],t_1<t_2$ such that
$\tilde{\varphi}(t_1)\in\partial(\mathcal{R})_{\delta_1},\tilde{\varphi}(t_2)\in\partial(\mathcal{R})_{\delta^*_2}$.
Therefore,  $S_{T}(\tilde{\varphi})\geq
S_{t_1t_2}(\tilde{\varphi})\geq V(\partial (\mathcal{R})_{\delta_1},\partial (\mathcal{R})_{\delta^*_2})
\geq s_0\geq  s_0\wedge s_1$. This proves the claim. As a result, $\{X^{\varepsilon}\in F_{2,z}\}\subset \{X^{\varepsilon}(0)=z,\rho_{T}(X^{\varepsilon},\mathbb{F}_{T}^{z}(0.9(s_0\wedge s_1)))\geq\delta^*_2-\delta_2\},\forall z\in \overline{(\mathcal{R})_{\delta_1}}\backslash (\mathcal{R})_{\delta_2}.$
By (\ref{ufldpupbb}) of ${\bf(II_u)}$,
there exists $\varepsilon_2>0$ such that for any $\varepsilon\in(0,\varepsilon_2)$
and $z\in \overline{(\mathcal{R})_{\delta_1}}\backslash (\mathcal{R})_{\delta_2}$, we have
\[
\mathbb{P}_z(X^{\varepsilon}_T\in (\mathcal{R})_{\delta_2})\leq
\mathbb{P}_z(X^{\varepsilon}\in F_{2,z})\leq
\mathbb{P}_{z}\{\rho_{T}(X^{\varepsilon},\mathbb{F}_{T}^{z}(0.9(s_0\wedge s_1)))\geq\delta^*_2-\delta_2\}
\leq \exp\{-0.8(s_0\wedge s_1)/\varepsilon^2\}.
\]
Therefore the inequality (\ref{repconout2}) follows immediately.

The proof of (\ref{repconout3}) remains unchanged, as one may check below. Let $C=\{\varphi\in {\bf C}_T: \varphi(T)\notin (\mathcal{R})_{\delta_2}\},
G=\{\varphi\in {\bf C}_T: \rho_T(\varphi,\psi^{x})<\delta_*, x\in \overline{(\mathcal{R})_{\delta_2}}\}$,
where $\delta_*=\delta-\delta_2$. Then $G$ is an open set in ${\bf C}_T$
and $G\subset C$ by (\ref{reploc}).
Furthermore, $\sup_{z\in \overline{(\mathcal{R})_{\delta_2}}}S_{T}^{z}(G)\leq
\sup_{z\in \overline{(\mathcal{R})_{\delta_2}}}S_{T}(\psi^{z})\leq 0.1(s_0\wedge s_1)$ by (\ref{repup}).
By ${\bf(I_u)}$,  there exists $\varepsilon_3>0$ such
that for any $\varepsilon\in(0,\varepsilon_3)$ and $z\in\overline{(\mathcal{R})_{\delta_2}}$,
we have
\begin{align*}
  \mathbb{P}_z(X^{\varepsilon}_T\notin (\mathcal{R})_{\delta_2})
=&\mathbb{P}_z(X^{\varepsilon}\in C)\\
\geq& \mathbb{P}_z(X^{\varepsilon}\in G)\\
\geq&\exp\big\{-\big(\sup_{z\in \overline{(\mathcal{R})_{\delta_2}}}S_T(\psi^z)+0.1(s_0\wedge s_1)\big)/\varepsilon^2\big\}\\
\geq&\exp\{-0.2(s_0\wedge s_1)/\varepsilon^2\}.
\end{align*}
where we have used ${\bf(I_u)}$ in the second inequality $\ge$.
 From the above obtained inequality, we get that for any $\varepsilon\in(0,\varepsilon_3)$ and $z\in\overline{(\mathcal{R})_{\delta_2}}$,
\begin{align*}
 \mathbb{P}_z(X^{\varepsilon}_T\in (\mathcal{R})_{\delta_2})
=&1- \mathbb{P}_z(X^{\varepsilon}_T\notin (\mathcal{R})_{\delta_2})\\
\leq &1-\exp\{-0.2(s_0\wedge s_1)/\varepsilon^2\}.
\end{align*}
So (\ref{repconout3}) follows easily from the above fact. This completes the proof.

Here is the proof of Theorem \ref{mainatt} with FWULDP.

{\bf Proof of Theorem \ref{mainatt} Under FWULDP}.
Let $\delta'\in(0,{\rm dist}\left(y,\mathcal{A}\right)/2)$ such that
$(\mathcal{A})_{\delta'}$ is a fundamental neighborhood of $\mathcal{A}$.
By Lemma \ref{klofA}, there exist $\delta_1>\delta^*_1>\delta_2>\delta_3, \delta_1,\delta^*_1,\delta_2,\delta_3\in(0,\delta')$
and $s_0>0$ such that $V\left(\partial (\mathcal{A})_{\delta_2},\partial (\mathcal{A})_{\delta^*_1}\right)\geq s_0$.
Since $\mathcal{A}$ is an attractor, we obtain
$T_0:=\inf\big\{u\geq 0:\Psi_{t}\big(\overline{(\mathcal{A})_{\delta_1}}\big)\subset(\mathcal{A})_{\delta_3}, t\geq u\big\}<+\infty$.
The set $F_0=\{\varphi\in {\bf C}_{T_0}:\varphi(0)\in \overline{(\mathcal{A})_{\delta_1}},\varphi(T_0)\in ((\mathcal{A})_{\delta_2})^c\}$ is a closed subset of ${\bf C}_{T_0}$. $F_0^0$ is bounded and $F_0$ does not contain any solution of system (\ref{unpersys}).
Thus, by Hypothesis \ref{closeinfH}, $s_1:=S_{T_0}(F_0)>0$.
Let $\delta=\delta_2\wedge \frac{s_0\wedge s_1}{10L}$,
where $L=L_{\bar{B}_{\delta'}(y)\cup \overline{(\mathcal{A})_{\delta'}}}$ is a constant as Lemma \ref{Vctin}.
Since $V(y,\mathcal{A})=0$, there exist $T_1>0$ and $ \tilde{\varphi}\in {\bf C}_{T_1}$ with $
\tilde{\varphi}(0)=y$ and $\tilde{\varphi}(T_1)\in \mathcal{A}$ such that $S_{T_1}(\tilde{\varphi})<0.2(s_0\wedge s_1)$.

As in the proof of Theorem \ref{mainatt} Under DZULDP, we only have prove (\ref{saddleep10}) and (\ref{saddleep20}). Here we replace $0.9$ in (\ref{saddleep20}) by $0.8$ and $\kappa$ in (\ref{Expdec})  by $0.4(s_0\wedge s_1)$.

To prove (\ref{saddleep10}), let
$G^z=\{\varphi\in {\bf C}^z_{T}:\rho_T(\varphi,\psi^z)<\delta_1\}$ for every $z\in \bar{B}_{\delta}(y)$. Fix $z\in \bar{B}_{\delta}(y)$. Recall $\psi^z(T)\in \mathcal{A}$, we get $\{X^{\varepsilon}_0=z,X^{\varepsilon}_T\in(\mathcal{A})_{\delta_1}\}\supset \{X^{\varepsilon}_0=z,X^{\varepsilon}\in G^z\}$.
By ${\bf(I_u)}$,  there exists $\varepsilon_1>0$ such that for any $\varepsilon\in(0,\varepsilon_1),z\in\bar{B}_{\delta}(y)$ we have
\begin{equation}\label{saddleep110_dis}
  \mathbb{P}_{z}\big(X^{\varepsilon}_T\in(\mathcal{A})_{\delta_1} \big)
\geq \mathbb{P}_{z}\big(X^{\varepsilon}\in G^z\big)\geq \exp\big\{-\big(0.3(s_0\wedge s_1)+0.1(s_0\wedge s_1)\big)/\varepsilon^2\big\}
=\exp\{-0.4(s_0\wedge s_1)/\varepsilon^2\}.
\end{equation}
Hence (\ref{saddleep10}) follows from (\ref{saddleep110_dis}).

For any $z\in \overline{(\mathcal{A})_{\delta_1}}$, let $F_z=\{\varphi\in\mathbf{C}^z_T:\varphi(T)\notin (\mathcal{A})_{\delta_1}\}$. Next we show that $\inf_{z\in\overline{(\mathcal{A})_{\delta_1}}}S^{z}_T(\{\varphi\in \mathbf{C}^z_T:\rho_T(\varphi,F_z)<\delta_1-\delta^*_1\})\geq s_0\wedge s_1$. In fact, fix $\tilde{z}\in \overline{(\mathcal{A})_{\delta_1}}$, let $\tilde{\varphi}\in\{\varphi\in \mathbf{C}^{\tilde{z}}_T:\rho_T(\varphi,F_{\tilde{z}})<\delta_1-\delta^*_1\}$. Because of the definition of $F_{\tilde{z}}$, we know that $\tilde{\varphi}(T)\notin (\mathcal{A})_{\delta^*_1}$. Either $S_{T}^{\tilde{z}}(\tilde{\varphi})\ge s_1$ or $S_{T}^{\tilde{z}}(\tilde{\varphi})< s_1$. Let us consider the later.
Then $S_{T_0}^{\tilde{z}}(\tilde{\varphi})<s_1$. So $\tilde{\varphi}(T_0)\in(\mathcal{A})_{\delta_2}$ from the definition of $s_1$.
Note that $\tilde{\varphi}(T)\notin (\mathcal{A})_{\delta^*_1}$
and $\overline{(\mathcal{A})_{\delta_2}}\subset (\mathcal{A})_{\delta^*_1}$.
Then by the continuity of $\tilde{\varphi}$, there exist $t_1,t_2\in(T_0,T],t_1>t_2$
such that $\tilde{\varphi}(t_1)\in \partial(\mathcal{A})_{\delta^*_1},\tilde{\varphi}(t_2)\in \partial(\mathcal{A})_{\delta_2}$.
Thus $S_{T}^{\tilde{z}}(\tilde{\varphi})\geq
S_{t_2t_1}^{\tilde{z}}(\tilde{\varphi})\geq
V(\partial (\mathcal{A})_{\delta_2},\partial (\mathcal{A})_{\delta^*_1})\geq s_0$.
Therefore $\inf_{z\in\overline{(\mathcal{A})_{\delta_1}}}S^{z}_T(\{\varphi\in \mathbf{C}^z_T:\rho_T(\varphi,F_z)<\delta_1-\delta^*_1\})\geq s_0\wedge s_1$. As a result, $\{X^{\varepsilon}\in F_z\}\subset \{X^{\varepsilon}(0)=z,\rho_{T}(X^{\varepsilon},\mathbb{F}_{T}^{z}(0.9(s_0\wedge s_1)))\geq\delta_1-\delta^*_1\},\forall z\in \overline{(\mathcal{A})_{\delta_1}}.$
Finally, by (\ref{ufldpupbb}) of ${\bf(II_u)}$,
there exists $\varepsilon_2>0$ such that for any
$\varepsilon\in(0,\varepsilon_2),z\in\overline{(\mathcal{A})_{\delta_1}}$ we have
\[\mathbb{P}_z\big(X^\varepsilon_T\notin (\mathcal{A})_{\delta_1}\big)
=\mathbb{P}_z\big(X^\varepsilon\in F_z\big)\leq
\mathbb{P}_{z}\{\rho_{T}(X^{\varepsilon},\mathbb{F}_{T}^{z}(0.9(s_0\wedge s_1)))\geq\delta_1-\delta^*_1\}
\leq\exp\{-0.8(s_0\wedge s_1)/\varepsilon^2\} . \]
This implies
\begin{align*}
   \int_{(\mathcal{A})_{\delta_1}}\mu^{\varepsilon}(dz)\mathbb{P}_z\big(X^\varepsilon_T\in (\mathcal{A})_{\delta_1}\big)=&
 \int_{(\mathcal{A})_{\delta_1}}\mu^{\varepsilon}(dz)\big(1-\mathbb{P}_z\big(X^\varepsilon_T\notin (\mathcal{A})_{\delta_1}\big)\big)\\
=&\mu^{\varepsilon}((\mathcal{A})_{\delta_1})
-\int_{(\mathcal{A})_{\delta_1}}\mu^{\varepsilon}(dz)\mathbb{P}_z\big(X^\varepsilon_T\notin (\mathcal{A})_{\delta_1}\big)\\
\geq& \mu^{\varepsilon}((\mathcal{A})_{\delta_1}) -\mu^{\varepsilon}((\mathcal{A})_{\delta_1})
\exp\{-0.8(s_0\wedge s_1)/\varepsilon^2\}
\end{align*}
for any $\varepsilon\in(0,\varepsilon_2)$, i.e., the inequality (\ref{saddleep20}) follows.

\vskip 0.2cm
\noindent{\bf  Acknowledgements}\quad  This work is supported by the National Natural Science Foundation of China (Nos. 12171321, 12001373, 11771295, 11931004).

%\def\refname{References}

%\bibliographystyle{imsart-number}
% Style BST file (imsart-number.bst or imsart-nameyear.bst)
\bibliography{main}

\end{document}